\documentclass[11pt]{amsart}
\usepackage{amsfonts,amssymb, enumerate,youngtab}

\input prepictex
\input pictex
\input postpictex

\numberwithin{equation}{section}

\newtheorem{theorem}{Theorem}[section]
\newtheorem{lemma}[theorem]{Lemma}

\newtheorem{corollary}[theorem]{Corollary}

\theoremstyle{definition}

\theoremstyle{remark}

\newcommand{\CC}{\mathbb{C}}
\newcommand{\HH}{\mathbb{H}}
\newcommand{\RR}{\mathbb{R}}

\newcommand{\ZZ}{\mathbb{Z}}

\newcommand{\OO}{\mathcal{O}}

\newcommand{\hh}{\mathfrak{h}}

\newcommand{\la}{\langle}
\newcommand{\ra}{\rangle}

\newcommand{\ttt}{\mathfrak{t}}

\begin{document}
\title[Jack polynomials and Cherednik algebras]{Generalized Jack
polynomials and the
representation theory of rational Cherednik algebras}
\author{Charles Dunkl}

\address{Department of Mathematics\\  
University of Virginia\\
P.O. Box 400137\\
Charlottesville, VA 22904-4137\\
}
\email{cfd5z@virginia.edu}

\author{Stephen Griffeth}
\address{School of Mathematics\\
University of Edinburgh\\
Edinburgh EH9 3JZ}
\email{griffeth@ed.ac.uk}

\begin{abstract}
We apply the Dunkl-Opdam operators and generalized Jack
polynomials to study category $\OO_c$ for
the rational Cherednik algebra of type $G(r,p,n)$.  We determine the
set of aspherical values and, in case $p=1$, answer a question of Iain Gordon on the ordering of category $\OO_c$.  
\end{abstract}

\thanks{We are very grateful to Pavel Etingof for
pointing out Lemma~\ref{aspherical lemma}, which helped to put us on
the right track, to an anonymous referee for careful reading and helpful advice, and to Iain Gordon for stimulating conversation and useful input.  The second
author acknowledges the full financial support of EPSRC grant EP/G007632.}

\maketitle

\section{Introduction}

We study the rational Cherednik algebra $\HH_c$
of type $G(r,1,n)$ by means of the Dunkl-Opdam operators introduced in
\cite{DuOp} and the generalized Jack polynomials
introduced in \cite{Gri}.  Our main results are a characterization of
the set of parameters $c$ for which $\HH_c e_\xi \HH_c \neq
\HH_c$, where $\xi$ is a linear character of $G(r,1,n)$ and $e_\xi$ is the corresponding idempotent, and a description, in terms of partition combinatorics, of an order on category $\OO_c$ that arises by studying the generalized eigenspaces of the Dunkl-Opdam operators on standard modules.  

Parameters $c$ such that $\HH_c e_\xi \HH_c \neq \HH_c$ are called $\xi$-\emph{aspherical}, and it is a standard fact that the functor $M \mapsto e_\xi M$ gives an equivalence from $\HH_c \mathrm{-mod}$ to $e_\xi \HH_c e_\xi \mathrm{-mod}$ exactly if $c$ is \emph{not} $\xi$-aspherical.  Writing $e$ for the symmetrizing idempotent (corresponding to the trivial character), Proposition~5.5 of Berest-Chalykh \cite{BeCh} states that for certain choices of $\xi$, we have an isomorphism $e \HH_c e \cong e_\xi \HH_{c'} e_\xi$ for a parameter $c'$ depending on $c$ and $\xi$, so combining this with our results allows one to establish Morita equivalences between rational Cherednik algebras at different parameters.  We hope that by combining these equivalences with the techniques of the papers \cite{Gor}, \cite{GoSt} and \cite{GGS} one may obtain a tighter relationship between the Cherednik algebra of type $G(r,1,n)$ and certain quiver varieties, including the Hilbert schemes of points on resolutions of type $A$ surface singularities.

Our determination of the set of $\xi$-aspherical parameters is similar to that in \cite{Dun} for the case $r=1$; in this case the result was first proved by Gordon-Stafford \cite{GoSt} and
Bezrukavnikov-Etingof \cite{BeEt}.  The parameter $c$ is a tuple $c=(c_0,d_0,\dots,d_{r-1})$ of complex numbers with $d_0+d_1+\cdots+d_{r-1}=0$ (see \eqref{RCA def} for its relationship to $\HH_c$).  For any $l \in \ZZ$ we define $d_l=d_{l'}$ if $l=l'\ \mathrm{mod} \ r$ and $0 \leq l' \leq
r-1$.  When $\xi=\mathrm{triv}$ is the trivial character, we abbreviate $\xi$-aspherical values to \emph{aspherical values}.  Our first main result, Theorem~\ref{aspherical theorem}, implies that the set of aspherical values is a certain finite union of hyperplanes.  Here we restate it with the minimum amount of notation.
\begin{theorem} 
The parameter $c$ is aspherical exactly if 
\begin{enumerate}
\item[(a)] $c_0=-k/m$ for integers $k$
and $m$ satisfying $1 \leq k < m \leq n$ or
\item[(b)] there is an integer $0 \leq l \leq r-1$, an integer $-(n-1)
\leq m \leq n-1$, and an integer $k$ such that $k \neq 0\ \mathrm{mod} \ r$,
$$k=d_l-d_{l-k}+r m c_0, \quad \mathrm{and} \quad 1 \leq k \leq  l+\left(\sqrt{n+\frac{1}{4} m^2}-\frac{1}{2} m-1\right)r
$$
\end{enumerate}
\end{theorem}  As in Lemma~\ref{aspherical lemma}, the hyperplanes of part (b) may also be
described as follows: for each integer $0 \leq l \leq r-1$, each
rectangular partition $\lambda$ with at most $n$ boxes, and each
integer $k \neq 0 \ \mathrm{mod} \ r$ with $$1 \leq k \leq l+r
(\mathrm{length}(\lambda)-1)$$ the hyperplane 
$$k=d_l-d_{l-k}+r(\lambda_1-\mathrm{length}(\lambda)) c_0 $$ is aspherical.  See Corollary~\ref{twists} for other
linear characters $\xi$ of $G(r,1,n)$, and Corollary~\ref{p corollary} for the case of $G(r,p,n)$.   

The irreducible representations of $G(r,1,n)$ are in sensible bijection with $r$-partitions $\lambda^\bullet=(\lambda^0,\dots,\lambda^{r-1})$ of $n$, and attached to each irreducible representation $S^{\lambda^\bullet}$ of $G(r,1,n)$ is a \emph{standard module} $M_c(\lambda^\bullet)$ of $\HH_c$.  The proof of Theorem~\ref{aspherical theorem} consists of computing the norm, with
respect to the contravariant form, of a certain invariant polynomial
in each standard module.  This gives a necessary condition for $c$ to
be aspherical.  On the other hand we can show that all such
parameters are aspherical using Theorem~7.5 of \cite{Gri} and the
results of \cite{BeEt}.  

By using the calculation of the spectrum of standard modules for $\HH_c$ in Theorem~5.1 of \cite{Gri}, we also answer question 10.1 from
\cite{Gor}.  Gordon's question has to do with \emph{category} $\OO_c$ for the Cherednik algebra: this is the Serre subcategory of $\HH_c$-modules generated by the standard modules $M_c(\lambda^\bullet)$.  Category $\OO_c$ is a highest weight category with standard objects the modules $M_c(\lambda^\bullet)$, where the order comes from a certain grading element of $\HH_c$ (the \emph{deformed Euler field}).   

When $d_l/rc_0 \in \ZZ$ for $0 \leq l \leq r-1$ the tuple $$(d_{r-1}/rc_0,\dots,d_0/rc_0)$$ is an element of the root lattice of type $A_{r-1}$ and hence corresponds, according to the (standard) bijection described in \ref{cores subsection}, to an $r$-core
partition $\mu_c$.  Via the bijection between ordinary partitions with
$r$-core $\mu_c$ and $r$-partitions obtained by taking $r$-quotients, dominance order on ordinary partitions induces an ordering $\geq_c'$ on
$r$-partitions .  This ordering arose
geometrically from quiver varieties in \cite{Gor}; here we show that
it arises also from the spectrum of the Dunkl-Opdam subalgebra.  It should be quite interesting to investigate a conceptual reason for this coincidence.  

For any parameter $c=(c_0,d_0,\dots,d_{r-1}) \in \RR^{r+1}$ with $c_0>0$, we define, as in \eqref{dominance condition}, another ordering on $r$-partitions $\lambda^\bullet=(\lambda^0,\dots,\lambda^{r-1})$ by $\lambda^\bullet \geq_c \chi^\bullet$ if for all $j \in \RR$ and $0 \leq l \leq r-1$,
\begin{align*} 
&|\{b \in \lambda^\bullet \ | \ \frac{d_{\beta(b)}}{r c_0}+\mathrm{ct}(b) > j \ \mathrm{or} \ \frac{d_{\beta(b)}}{r c_0}+\mathrm{ct}(b)=j \ \mathrm{and} \ \beta(b) \leq l \}| \\
& \geq |\{b' \in \chi^\bullet \ | \ \frac{d_{\beta(b')}}{r
c_0}+\mathrm{ct}(b') > j \ \mathrm{or} \ \frac{d_{\beta(b')}}{r
c_0}+\mathrm{ct}(b')=j \ \mathrm{and} \ \beta(b') \leq l \}| \notag.
\end{align*}  For the definitions of the statistics $\mathrm{ct}(b)$ and $\beta(b)$ see \eqref{content def} and \eqref{beta def}.  With these notations, Theorem~\ref{ordering theorem} and Corollary~\ref{iainquestion} state that
\begin{theorem}
\begin{enumerate}
\item[(a)] For any choice of parameter $c \in \RR^{r+1}$ with $c_0 >0$, category $\OO_c$ is a highest weight category with respect to the ordering $\geq_c$ described above;
\item[(b)] if $d_l/r c_0 \in \ZZ$ for $0 \leq l \leq r-1$ then $\OO_c$ is a highest weight category with respect to the ordering $\geq_c'$ described above.  
\end{enumerate}
\end{theorem}  Since there are fewer order relations for these orders than for the order
coming from the deformed Euler field, our result broadens the applicability
of Theorem~4.49 of \cite{Rou}, allowing one to deduce Morita equivalences between category $\OO$'s at different parameters.

\subsection{The symmetric group}

Let $S_n$ be the group of permutations of the set $\{1,2,\dots,n\}$.  The notation $w_1 \leq w_2$ for $w_1,w_2 \in S_n$ refers to Bruhat order, and we write $l(w)$ for the length of an element $w \in S_n$.  Let $w_0 \in S_n$ be the \emph{longest element}, with
\begin{equation}
w_0(i)=n-i+1 \quad \hbox{for $1 \leq i \leq n$.}
\end{equation}  For a sequence $\mu \in \ZZ_{\geq 0}^n$ of $n$ non-negative integers, we write $\mu^+$ for the non-increasing (partition) rearrangement of $\mu$, and $\mu^-$ for the non-decreasing (anti-partition) rearrangement of $\mu$.  For $w \in S_n$ and $\mu \in \ZZ_{\geq 0}^n$, the formula
\begin{equation}
w.\mu=(\mu_{w^{-1}(1)},\mu_{w^{-1}(2)},\dots,\mu_{w^{-1}(n)})
\end{equation} defines a left action of $S_n$ on $\ZZ_{\geq 0}^n$.  Let $w_\mu$ be the longest element of $S_n$ such that $w_\mu.\mu=\mu^-$; thus
\begin{equation} \label{wmu def}
w_\mu(i)=|\{1 \leq j < i \ | \ \mu_j < \mu_i \}|+|\{i \leq j \leq n \ | \ \mu_j \leq \mu_i \}|.
\end{equation}  Also define the \emph{rank function} $r_\mu$ by
\begin{equation}
r_\mu(i)=|\{1 \leq j \leq i \ | \ \mu_j \geq \mu_i \}|+|\{i<j\leq n \ | \ \mu_j > \mu_i \}|
\end{equation} so that
\begin{equation}
w_\mu(i)+r_\mu(i)=n+1 \quad \text{or equivalently} \quad w_\mu=w_0 r_\mu.
\end{equation}  There is a partial order $\geq$ on $\ZZ_{\geq 0}^n$ defined as follows: for $\mu, \nu \in \ZZ_{\geq 0}^n$ we have
\begin{equation} \label{order def}
\mu > \nu \quad \text{if} \quad \mu_+ >_d \nu_+ \quad \text{or} \quad \mu_+=\nu_+ \quad \text{and} \quad w_\mu > w_\nu,
\end{equation} where $>_d$ is dominance order on partitions.  For a reminder of the definition of dominance order on partitions, see \eqref{dom def}; here we use the notation $>_d$ to distinguish it from the order on $\ZZ_{\geq 0}^n$ that we are defining.

\subsection{Partitions, Young diagrams, and tableaux}

Let $n \in \ZZ_{>0}$ be a positive integer.  A \emph{partition of
length $n$} is a non-increasing sequence $\lambda=(\lambda_1 \geq
\lambda_2 \geq \cdots \geq \lambda_n > 0)$ of $n$ positive integers.
Without reference to length, a \emph{partition} is a non-increasing
sequence $\lambda=(\lambda_1 \geq \lambda_2 \geq \cdots \geq \lambda_m
\geq 0)$ of non-negative integers, and we consider two partitions to
be equal if they differ by a terminal string of zeros.  In particular,
we refer to a sequence of $0$'s as the \emph{empty} partition
$\emptyset$.  Let $r$ be a positive integer.  An \emph{$r$-partition}
is a sequence
$\lambda^\bullet=(\lambda^0,\lambda^1,\dots,\lambda^{r-1})$ (some
$\lambda^i$'s may be empty) of $r$ partitions.  The \emph{Young
diagram} of an $r$-partition is the graphical representation consisting of a collection of boxes stacked in a corner: the Young diagram for the $4$-partition $\left((3,3,1),(2,1),\emptyset,(5,5,2,1) \right)$ of $23$ is 

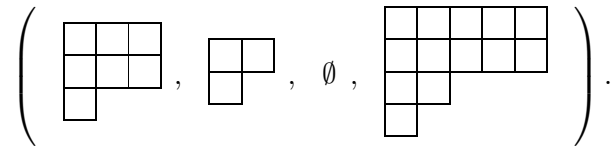
\begin{equation*}
\left(\begin{array}{cccc} \begin{array}{c} \yng(3,3,1) \end{array}, &  \begin{array}{@{}c} \yng(2,1) \end{array}, &  \begin{array}{@{}c} \emptyset \end{array}, &  \begin{array}{@{}c} \yng(5,5,2,1) \end{array}
\end{array} \right).
\end{equation*}   

A \emph{tableau} $T$ on an $r$-partition $\lambda^\bullet$ of $n$ is a filling
of the boxes of $\lambda^\bullet$ with non-negative integers.  A tableau $T$ is
\emph{column-strict} if within each partition $\lambda^i$, its entries
are strictly increasing from top to bottom and weakly increasing from
left to right.  A \emph{standard Young tableau} on $\lambda^\bullet$ is a
bijection $T$ from the boxes of $\lambda^\bullet$ to the set $\{1,2,\dots,n\}$ in
such a way that the entries are increasing left to right and top to bottom.  For example,
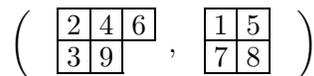
\begin{equation*}
\left( \begin{array}{cccc} \begin{array}{c} \young(246,39) \end{array},\ \begin{array}{c} \young(15,78) \end{array} \end{array} \right)
\end{equation*} is a standard Young tableau on the $2$-partition $((3,2),(2,2))$.   We write
\begin{equation}
T(b)=i \quad \hbox{if $i$ appears in box $b$.}
\end{equation}  Then for standard Young tableaux $S$ and $T$ on $\lambda$, the permutation $T S^{-1}$ of $\{1,2,\dots,n\}$ is a measure of the distance between $S$ and $T$.  For an $r$-partition $\lambda^\bullet$ we write
\begin{equation}
\mathrm{SYT}(\lambda)=\{\hbox{standard Young tableaux on $\lambda^\bullet$} \}
\end{equation}  

For a box $b \in \lambda$, write $\mathrm{row}(b)=i$ and $\mathrm{col}(b)=j$ if $b$ is in the $i$th row and $j$th column of $\lambda$, and define the \emph{content} $\mathrm{ct}(b)$ of $b$ by 
\begin{equation} \label{content def}
\mathrm{ct}(b)=\mathrm{col}(b)-\mathrm{row}(b).
\end{equation}  We also define the function $\beta$ on the set of boxes of $\lambda^\bullet$ by 
\begin{equation} \label{beta def}
\beta(b)=l \quad \hbox{if $b \in \lambda^l$.}
\end{equation}  Thus for the tableau $T$ pictured above, one has
\begin{equation*}
\beta(T^{-1}(7))=1 \quad \mathrm{and} \quad \mathrm{ct}(T^{-1}(3))=-1.
\end{equation*}

\subsection{The groups $G(r,1,n)$}

Fix positive integers $r$ and $n$.  A \emph{monomial matrix} is a square matrix with exactly one non-zero entry in each row and each column.  Let
\begin{equation}
W=G(r,1,n)
\end{equation} be the group of $n$ by $n$ monomial matrices whose non-zero entries are $r$th roots of $1$.  We write $\CC W$ for the complex group algebra of $W$.  Let $\zeta=e^{2 \pi i/r}$ and
\begin{equation}
\zeta_i=\mathrm{diag}(1,\dots,\zeta,\dots,1)
\end{equation} be the diagonal matrix with a $\zeta$ in the $i$th position and $1$'s elsewhere on the diagonal.  Let
\begin{equation}
s_{ij}=(ij) \quad \mathrm{and} \quad s_i=s_{i,i+1}
\end{equation} be the transposition interchanging $i$ and $j$ and the
$i$th simple transposition, respectively.  We write $\hh=\CC^n$ for
the defining representation of $W$ and $\hh^*$ for its dual.

There is a version of Young's orthonormal form that works for the
groups $G(r,1,n)$: the irreducible complex representations of
$G(r,1,n)$ may be indexed by the $r$-partitions of $n$ in such a way
that if $S^{\lambda^\bullet}$ is the irreducible representation
corresponding to the $r$-partition $\lambda^\bullet$, then
$S^{\lambda^\bullet}$ has a basis $v_T$ indexed by the set of standard
Young tableaux $T$ on $\lambda^\bullet$.  The vectors $v_T$ are
eigenvectors for a certain maximal commutative subalgebra of $\CC W$,
and the action of the generators $\zeta_i$ and $s_i$ on $v_T$ can be
made quite explicit: see section 3 of \cite{Gri} for these facts.

\subsection{The rational Cherednik algebra} \label{RCA def}

Let $c_0$ and $d_i$ for $i \in \ZZ$ be variables such that $d_i=d_j$
if $i=j\ \mathrm{mod} \ r$, and $d_0+d_1+\cdots+d_{r-1}=0$.  Let
$k=\CC[c_0,d_i]_{1 \leq i \leq r-1}$ be the polynomial ring over
$\CC$ generated by these variables and let $F$ be its fraction field.  The \emph{rational Cherednik algebra} $\HH$ for $W=G(r,1,n)$ is generated by $k[x_1,\dots,x_n]$, $k W$, and $k[y_1,\dots,y_n]$ with relations
\begin{equation}
y_i x_i=x_i y_i+1-c_0 \sum_{j \neq i} \sum_{l=0}^{r-1} \zeta_i^l s_{ij}\zeta_i^{-l}-\sum_{j=0}^{r-1} (d_j-d_{j-1}) e_{ij} 
\end{equation} and for $i \neq j$
\begin{equation}
y_i x_j=x_j y_i+c_0 \sum_{l=0}^{r-1} \zeta^{-l} \zeta_i^l s_{ij} \zeta_i^{-l},
\end{equation}  where for $1 \leq i \leq n$ and $0 \leq j \leq r-1$
\begin{equation}
e_{ij}=\frac{1}{r} \sum_{l=0}^{r-1} \zeta^{-lj} \zeta_i^l
\end{equation} are the idempotents for the cyclic reflection subgroups of $W$.  These parameters are related to those found in Gordon's paper \cite{Gor} by the equations
\begin{equation}
H_j=\frac{1}{r}(d_{j-1}-d_j) \quad \hbox{for $0 \leq j \leq r-1$, and} \quad h=-c_0,
\end{equation} and to those in section~6 of Rouquier's paper \cite{Rou} by
\begin{equation}
h_j=-\frac{1}{r} d_j \quad \hbox{for $0 \leq j \leq r-1$, and} \quad h=-c_0.
\end{equation}  Thus setting $q=e^{-2 \pi i c_0}$ and $Q_j=e^{-2 \pi i d_j/r}$ for $0 \leq j \leq r-1$ one obtains the parameters for the corresponding finite Hecke algebra of type $G(r,1,n)$.  Assuming $r>1$ and $n>1$ as in 6.1.2 of \cite{Rou}, this finite Hecke algebra has generators $T_0,T_1,\dots,T_{n-1}$, with the relations $T_0 T_1 T_0 T_1=T_1 T_0 T_1 T_0$, braid relations between $T_i$ and $T_j$ for $1 \leq i < j \leq n-1$, and 
$$(T_i-q)(T_i+1)=0 \quad \hbox{for $1 \leq i \leq n-1$, and} \quad \prod_{j=0}^{r-1} (T_0-e^{2 \pi i j/r} Q_j)=0.$$

The \emph{PBW theorem} for $\HH$ asserts that as a $k$-module
\begin{equation}
\HH \cong k[x_1,\dots,x_n] \otimes k W \otimes k[y_1,\dots,y_n].
\end{equation}  It implies that if we define the \emph{standard module} $M(\lambda^\bullet)$ by
\begin{equation}
M(\lambda^\bullet)=\mathrm{Ind}_{k[y_1,\dots,y_n] \rtimes kW}^\HH S^{\lambda^\bullet} 
\end{equation} then as a $k[x_1,\dots,x_n] \rtimes k W$-module,
\begin{equation}
M(\lambda^{\bullet}) \cong k[x_1,\dots,x_n] \otimes S^{\lambda^\bullet}.
\end{equation}  

Let $\overline{\cdot}$ be the automorphism of $k$ that fixes $c_0,d_0,\dots,d_{r-1}$ and acts as complex conjugation on $\CC$.  Let $x \mapsto \overline{x}$ be a conjugate-linear $W$-equivariant
isomorphism from $\hh^*$ onto $\hh$.  Fix a $W$-invariant positive definite Hermitian form on $S^{\lambda^\bullet}$.  Then $M(\lambda^\bullet)$ carries a \emph{contravariant form} $\la
\cdot,\cdot \ra$, determined uniquely by the requirements that it
be linear in the second variable, skew-symmetric, 
and satisfy the conditions
\begin{equation}
\la w.f,w.g \ra=\la f,g \ra \quad \mathrm{and} \quad \la x.f,g \ra=\la f,\overline{x}.g
\ra \quad \hbox{for $f,g \in M(\lambda^\bullet)$, $w \in W$, and $x
\in \hh^*$,}
\end{equation} and restrict to the given form on $S^{\lambda^\bullet}$.  The proof of Theorem 2.18 of \cite{DuOp} works, mutatis mutandis, for this situation.

If we choose a particular specialization $c_0,d_i \in \CC$ of the
parameters to complex numbers, we write $\HH_c$ for the resulting
algebra.  The specializations
$M_c(\lambda^\bullet)$ of the standard modules have unique irreducible
quotients $L_c(\lambda^\bullet)$, and we define \emph{category $\OO_c$} to be the
category of finitely generated $\HH_c$-modules on which
$y_1,\dots,y_n$ act locally nilpotently.  The contravariant form also specializes to a conjugate-symmetric form $\la \cdot, \cdot \ra_c$, provided that the parameters are real numbers.  Otherwise it is not conjugate-symmetric; see Theorem 2.18 of \cite{DuOp}.

\subsection{The Dunkl-Opdam subalgebra and the non-symmetric generalized Jack polynomials}

For $1 \leq i \leq n$ define
\begin{equation}
z_i=y_i x_i+c_0 \phi_i \ \mathrm{where} \ \phi_i=\sum_{1 \leq j <i} \zeta_i^l s_{ij} \zeta_i^{-l} 
\ \hbox{are the Jucys-Murphy elements for $G(r,1,n)$.}
\end{equation}   The elements $z_i$ were introduced in \cite{DuOp},
where it is also proved that they commute.  The \emph{Dunkl-Opdam} subalgebra $\ttt$ of $\HH$ is generated by $z_1,\dots,z_n$ and $\zeta_1,\dots,\zeta_n$.  By the PBW theorem it is isomorphic to the polynomial ring in the variables $z_1,\dots,z_n$ tensored with the group algebra of $(\ZZ/r \ZZ)^n$.  Since the operators $z_i$ are fixed by the conjugate-linear anti-involution $\overline{\cdot}$ of $\HH$ such that $\overline{x_i}=y_i$ and $\overline{w}=w^{-1}$ for all $1 \leq i \leq n$ and $w \in G(r,1,n)$, they act as self-adjoint operators on $M(\lambda^\bullet)$ with respect to the contravariant form.  

The following theorem combines Theorems 5.1 and 6.1 of
\cite{Gri}.  The $T^{-1}$'s appear here because we are using the
``inverse'' definition of standard Young tableau to that in
\cite{Gri}.  Define a partial order on pairs $(\mu,T) \in \ZZ_{\geq 0}^n \times \mathrm{SYT}(\lambda^\bullet)$ by $(\mu,T) > (\nu,S)$ exactly if $\mu > \nu$, where the order on $\ZZ_{\geq 0}^n$ is defined by \eqref{order def}.  

\begin{theorem} \label{Upper triangular}
Let $\lambda^\bullet$ be an $r$-partition of $n$, $\mu \in \ZZ_{\geq 0}^n$, and let $T$ be a standard tableau on $\lambda^\bullet$.  Put $v_T^\mu=w_\mu^{-1}.v_T$ and recall the definitions of $\beta$ and ct given in \eqref{content def} and \eqref{beta def}.
\begin{enumerate}
\item[(a)]  The action of $\zeta_i$ and $z_i$ on $M(\lambda^\bullet)$ are given by
\begin{align*}
\zeta_i.x^\mu v_T^\mu=\zeta^{\beta(T ^{-1}w_\mu(i))-\mu_i} x^\mu v_T^\mu
\end{align*}
and
\begin{align*} \label{z spectrum}
z_i.x^\mu v_T^\mu&=\left( \mu_i+1-(d_{\beta(T^{-1}
w_\mu(i))}-d_{\beta(T^{-1} w_\mu(i))-\mu_i-1})- r \mathrm{ct}(T^{-1}w_\mu(i))c_0 \right) x^\mu v_T^\mu \\
&+\sum_{(\nu,S) < (\mu,T)} c_{\nu,S} x^\nu v_{S}^\nu.
\end{align*}
\item[(b)]  Assuming that scalars are extended to $F=\CC(c_0,d_1,d_2,\dots,d_{r-1})$, for each $\mu \in \ZZ_{\geq 0}^n$ and $T \in \mathrm{SYT}(\lambda^\bullet)$ there exists a unique $\ttt$ eigenvector $f_{\mu,T} \in M(\lambda)$ such that
\begin{equation*}
f_{\mu,T}=x^\mu v_T^\mu+\mathrm{lower} \ \mathrm{terms}.
\end{equation*}  
The $\ttt$-eigenvalue of $f_{\mu,T}$ is determined by the formulas in
part (a).
\item[(c)] Write
\begin{equation*}
a_i=\mathrm{ct}(T^{-1}w_\mu(i)) \quad \mathrm{and} \quad b_i=\beta(T^{-1}w_\mu(i)).
\end{equation*}  Then the norm of $f_{\mu,T}$ is given by
\begin{align*}
\la f_{\mu,T}, f_{\mu,T} \ra&=\prod_{i=1}^n \prod_{k=1}^{\mu_i} \left( k-(d_{b_i}-d_{b_i-k})- r a_i c_0\right) \\
&\times \prod_{\substack{1 \leq i < j \leq n \\ \mu_i > \mu_j}} \prod_{\substack{1 \leq k \leq \mu_i-\mu_j \\ k=b_i-b_j \ \mathrm{mod} \ r}} \frac{\left(k-(d_{b_{i}}-d_{b_{j}})-r (a_{i}-a_{j})c_0\right)^2-(rc_0)^2}{\left(k-(d_{b_{i}}-d_{b_{j}})-r (a_{i}-a_{j})c_0 \right)^2  } \\
&\times  \prod_{\substack{1 \leq i < j \leq n \\ \mu_i < \mu_j-1}} \prod_{\substack{1 \leq k \leq \mu_j-\mu_i-1 \\ k=b_j-b_i \ \mathrm{mod} \ r}} \frac{\left(k-(d_{b_{j}}-d_{b_{i}})-r (a_{j}-a_{i}) c_0 \right)^2-(rc_0)^2}{\left(k-(d_{b_{j}}-d_{b_{i}})-r(a_{j}-a_{i})c_0 \right)^2  }.
\end{align*} 
\end{enumerate}
\end{theorem}

\section{Generalized Jack polynomials}

\subsection{Definition of the (symmetric) generalized Jack polynomials} The functions $f_{\mu,T}$ are a basis of $M(\lambda^\bullet)$.  Let
$e=\sum_{w \in S_n} w$ be the $S_n$-symmetrizer.  Then the functions
$e.f_{\mu,T}$ such that $\mu_i=\beta(T^{-1} w_\mu(i))\ \mathrm{mod} \ r$ for $1 \leq i \leq n$ span the $G(r,1,n)$ invariants in $M(\lambda^\bullet)$.  To single a basis out from this spanning set we use the following construction from \cite{Dun}.

Given $(\mu,T) \in \ZZ_{\geq 0}^n \times \mathrm{SYT}(\lambda^\bullet)$, define a function $S=S(\mu,T)$ on the boxes of $\lambda^\bullet$ by $S(b)=\mu_{w_\mu^{-1} T(b)}$.  It follows from the definition that $S$ is weakly increasing top to bottom and left to right.  For instance, if
\begin{equation*}
\mu=(2,3,2,0,4,2,5,2,2) \ \hbox{with $w_\mu=(6,7,5,1,8,4,9,3,2)$}
\end{equation*} and 
\begin{equation*}
T=\left( \begin{array}{cccc} \begin{array}{c} \young(134,89) \end{array}, \ \begin{array}{c}\young(26,57) \end{array} \end{array} \right) \qquad \mathrm{then} \qquad S(\mu,T)=\left( \begin{array}{cccc} \begin{array}{c} \young(022,45) \end{array},\ \begin{array}{c} \young(22,23) \end{array} \end{array} \right).
\end{equation*}  One checks that $\mu_i=\beta(T^{-1} w_\mu(i))$ mod
$r$ for all $1 \leq i \leq n$ exactly if $\beta(b)=S(b)\ \mathrm{mod} \ r$ for
all $b \in \lambda^\bullet$.  Therefore the invariants in
$M(\lambda^\bullet)$ for the
diagonal subgroup of $G(r,1,n)$ generated by $\zeta_1,\dots,\zeta_n$
are those $f_{\mu,T}$ so that $S(b)=\beta(b)\ \mathrm{mod} \ r$ for all $b \in
\lambda^\bullet$, where $S=S(\mu,T)$.  We will describe the $G(r,1,n)$
invariants.

The (Cherednik-style) intertwining operators are
\begin{equation}
\sigma_i=s_i+f_{\alpha_i}, \quad \text{where} \quad f_{\alpha_i}=\frac{c_0}{z_i-z_{i+1}} \sum_{l=0}^{r-1} \zeta_i^l \zeta_{i+1}^{-l}.
\end{equation}  Direct calculation shows
the $\sigma_i$'s satisfy the braid relations, and that
\begin{equation}  \label{sigma squared}
\sigma_i^2=1-f_{\alpha_i}^2
\end{equation} as in Lemma~5.2~(a) of \cite{Gri2}.  For a permutation $w \in S_n$ written as a reduced word in the simple reflections $w=s_{i_1} \cdots s_{i_p}$ we define
$\sigma_w=\sigma_{i_1} \cdots \sigma_{i_p}$; by the braid relations this does not depend on our choice of reduced word for $w$.  It follows from their definition that the intertwiners are self-adjoint for the contravariant form on $M(\lambda^\bullet)$.

\begin{lemma} \label{jacks lemma}
\begin{enumerate}
\item[(a)] If $S$ is a column-strict tableau on $\lambda^\bullet$ then there is a non-decreasing $\mu \in \ZZ_{\geq 0}^n$ and $T \in \mathrm{SYT}(\lambda^\bullet)$ with $S=S(\mu,T)$.
\item[(b)]  The $F$-span of the polynomials $e.f_{\mu,T}$ as $T$ ranges over $\mathrm{SYT}(\lambda^\bullet)$ and $\mu$ ranges over non-decreasing sequences is all of $M(\lambda^\bullet)^{S_n}$.
\item[(c)]  If $\mu \in \ZZ_{\geq 0}^n$ is non-decreasing and $T_1,T_2
\in \mathrm{SYT}(\lambda^\bullet)$ with $S(\mu,T_1)=S(\mu,T_2)$ then
$e.f_{\mu,T_1}$ and $e.f_{\mu,T_2}$ are $\CC$-multiples of one
another.
\item[(d)]  If $\mu \in \ZZ_{\geq 0}^n$ is non-decreasing, $T \in \mathrm{SYT}(\lambda^\bullet)$ and $S(\mu,T)$ is not column-strict then $e.f_{\mu,T}=0$. 
\end{enumerate}
\end{lemma}  
\begin{proof}
Part (a) follows from the definitions given above.  For part (b), the formula given in part (a) of Lemma~5.3 of \cite{Gri} implies that for each $(\nu,T) \in \ZZ_{\geq 0}^n \times \mathrm{SYT}(\lambda^\bullet)$ we have $f_{\nu,T}=\sigma_w f_{\nu^-,T}$ for some $w \in S_n$.  But $\sigma_w.f_{\nu^-,T} \in F G(r,1,n) f_{\nu^-,T}$ by definition of $\sigma_w$, so $e.f_{\nu,T}$ is an $F$-multiple of $e.f_{\nu^-,T}$, proving (b).    

For (c), assuming $S(\mu,T_1)=S(\mu,T_2)$ we obtain $R(\mu,T_1)=R(\mu,T_2)$ where $R(\mu,T)$ is defined in 7.17 of \cite{Gri}.  Hence Lemma 7.4 of \cite{Gri} shows that there is a sequence of intertwiners $\sigma_{i_1},\dots,\sigma_{i_l}$ with $f_{\mu,T_1}=\sigma_{i_l} \cdots \sigma_{i_1} f_{\mu,T_2}$ and so that $\sigma_{i_j} \cdots \sigma_{i_1} f_{\mu,T_2}=f_{\mu,S_j}$ for standard Young tableaux $S_j$ with the property that $R(\mu,S_j)=R(\mu,T_2)$ for all $1 \leq j \leq l$.  We claim that $f_{\mu,S_j} \in \CC S_n.f_{\mu,T_2}$ for all $j$; part (c) follows immediately from this.  We may therefore assume $l=1$, so that $f_{\mu,T_1}=\sigma_i.f_{\mu,T_2}$ and $\mu_i=\mu_{i+1}$.  But $\sigma_i=s_i+f_{\alpha_i}$, and the condition $\mu_i=\mu_{i+1}$ together with the formula in Theorem~\ref{Upper triangular} for the $\ttt$-eigenvalue of $f_{\mu,T_2}$ imply that $f_{\alpha_i}$ acts by a complex number on $f_{\mu,T_2}$, finishing the proof of (c).

We now prove (d).  Assuming $S=S(\mu,T)$ is not column-strict, there are boxes $b,b' \in \lambda^\bullet$ so that $b'$ is directly below $b$ and $S(b')=S(b)$.  By part (c) we may prove (d) for any tableau $T'$ with $S(\mu,T)=S(\mu,T')$.  We may therefore assume that $w_\mu^{-1} T(b)=i+1$ and $w_\mu^{-1} T(b')=i$ for some integer $1 \leq i \leq n-1$.  Using part (c) of Lemma~5.3 of \cite{Gri} implies $\sigma_i.f_{\mu,T}=0$, and the formula in Theorem~\ref{Upper triangular} for the $\ttt$-weight of $f_{\mu,T}$ implies that $s_i.f_{\mu,T}=-f_{\mu,T}$ whence $e.f_{\mu,T}=es_i.f_{\mu,T}=-e.f_{\mu,T}$.  Thus $e.f_{\mu,T}=0$.
\end{proof}
For each column-strict tableau $S$ on $\lambda^\bullet$ with $S(b)=\beta(b)\ \mathrm{mod} \ r$ let $\mu \in \ZZ_{\geq 0}^n$ be the sequence obtained by arranging the entries of $S$ in non-decreasing order, and fix a tableau $T \in \mathrm{SYT}(\lambda^\bullet)$ with $S=S(\mu,T)$.  Define the \emph{generalized Jack polynomial} by
\begin{equation}
g_S=e.f_{\mu,T}.
\end{equation}  The definition of $g_S$ depends on the choice of $T$,
but only up to multiplication by a constant in $\CC$.

\subsection{The norm of $g_S$}
In the following theorem, and all that follows, we use the convention
that the product over an empty set is $1$.  
\begin{theorem} \label{norm formula}
The functions $g_S$, where $S$ ranges over column strict tableau
satisfying $S(b)=\beta(b)\ \mathrm{mod} \ r$ for all $b \in \lambda^\bullet$,
are a basis of $M(\lambda^\bullet)^W$.  The norm of $g_S$ is given by
\begin{align*}
\la g_S,g_S \ra&=n! \prod_{b \in \lambda^\bullet} \prod_{1 \leq k \leq
S(b)} \left(k-(d_{\beta(b)}-d_{\beta(b)-k})-r \mathrm{ct}(b) c_0 \right)
\\
&\times \prod_{b,b' \in \lambda^\bullet} \prod_{\substack{1 \leq k \leq
S(b)-S(b') \\ k=\beta(b)-\beta(b') \ \mathrm{mod} \ r }}
\frac{k-(d_{\beta(b)}-d_{\beta(b')})-r(\mathrm{ct}(b)-\mathrm{ct}(b')-1) c_0
}{k-(d_{\beta(b)}-d_{\beta(b')})-r(\mathrm{ct}(b)-\mathrm{ct}(b'))c_0} \\
&\times \prod_{b, b' \in \lambda^\bullet} \prod_{\substack{1 \leq k \leq
S(b)-S(b')-r \\ k=\beta(b)-\beta(b') \ \mathrm{mod} \ r }}
\frac{k-(d_{\beta(b)}-d_{\beta(b')})-r(\mathrm{ct}(b)-\mathrm{ct}(b')+1) c_0
}{k-(d_{\beta(b)}-d_{\beta(b')})-r(\mathrm{ct}(b)-\mathrm{ct}(b'))c_0}
\end{align*}
\end{theorem}
\begin{proof}
We first show that the functions $g_S$ have distinct $\ttt^{S_n}$ eigenvalues; the fact that they  are linearly independent follows from this.  The eigenvalues of symmetric functions in $z_1,\dots,z_n$ acting on $g_S$ are determined by the multiset of eigenvalues of $z_1,\dots,z_n$ acting on any $f_{\mu,T}$ for which $S=S(\mu,T)$.  By the definition of $S(\mu,T)$ given above and part (a) of Theorem~\ref{Upper triangular} this multiset is
$$\{ \mu_i+1-(d_{\beta(T^{-1}w_\mu(i))}-d_{-1})-r \mathrm{ct}(T^{-1} w_\mu(i))c_0 \ | \ 1 \leq i \leq n\}=\{S(b)+1-(d_{\beta(b)}-d_{-1})-r\mathrm{ct}(b)c_0 \ | \ b \in \lambda^\bullet \}$$ and $S$ is uniquely determined by it, because any filling of a partition which is weakly increasing left to right and top to bottom is uniquely determined by the associated multiset of contents.

We next calculate the norm of $g_S=e.f_{\mu,T}$ using a version of the
argument proving (5.8.24) in \cite{Mac2}.  We write
$R=\{\epsilon_i-\epsilon_j \ | \ 1 \leq i \neq j \leq n \}$ for the
root system of type $A_{n-1}$ and $R^+=\{\epsilon_i-\epsilon_j \ | \ 1
\leq i < j \leq n \}$ for the set of positive roots.  For a root $\alpha=\epsilon_i-\epsilon_j$ of $S_n$ write
\begin{equation}
f_\alpha=\frac{c_0}{z_i-z_j} \sum_{l=0}^{r-1} (\zeta_i \zeta_j^{-1})^l
\end{equation} and let $g_\alpha$ be the constant by which $f_\alpha$
acts on $f_{\mu,T}$.  This is well-defined by the column-strictness of
$S$: setting $b_1=T^{-1} w_\mu(i)$ and $b_2=T^{-1} w_\mu(j)$, $f_\alpha$ acts by $0$ on $f_{\mu,T}$ unless
$\zeta^{\beta(b_1)-\mu_i}=\zeta^{\beta(b_2)-\mu_j}$, and
then the action is well-defined unless also 
$$\mu_i-\mu_j=d_{\beta(b_1)}-d_{\beta(b_2)}+r
(\mathrm{ct}(b_1)-\mathrm{ct}(b_2))c_0.$$  But then
$\beta(b_1)=\beta(b_2)$ and $\mathrm{ct}(b_1)=\mathrm{ct}(b_2)$ implies
$0=\mu_i-\mu_j=S(b_1)-S(b_2)$, contradicting the column-strictness of
$S$.   

We write $e=\sum_{w \in S_n} w$ and
\begin{equation}
e=\sum_{w \in S_n} \sigma_w f_w
\end{equation} for some $f_w$ in the ring obtained from $\ttt$ by inverting $z_i-z_j$ for $1 \leq i<j \leq n$.  Our first goal is to compute
each $f_w$ explicitly.  Observe that $f_{w_0}=1$ since no $\sigma_w$ other
than $\sigma_{w_0}$ contributes $w_0$.  Next, for each $1 \leq i \leq
n-1$ we have $s_i e=e$ and hence 
\begin{align*}
\sum_{w \in S_n} \sigma_w f_w
&=\left(\sigma_i-f_{\alpha_i} \right)\sum_{w \in S_n} \sigma_w f_w
=\sum_{s_i w>w} \sigma_{s_i w} f_w+\sum_{s_i w<w} \sigma_i^2
\sigma_{s_i w} f_w-\sum_{w \in S_n} \sigma_w f_{w^{-1} \alpha_i} f_w \\
&=\sum_{s_i w<w} \sigma_w f_{s_i w}+\sum_{s_i w>w}
\sigma_{w} f_w'-\sum_{w \in S_n} \sigma_w f_{w^{-1}\alpha_i} f_w
\end{align*} for some $f_w' \in \ttt$, where in the first and second summations we have reindexed, replacing $w$ by $s_i w$.  Thus for $w \in S_n$ with $s_i
w< w$ comparing coefficients of $\sigma_w$ gives the recurrence
\begin{equation}
f_w=f_{s_iw}-f_{w^{-1} \alpha_i} f_w \quad \implies \quad 
f_{s_i w}=\left(1-f_{-w^{-1} \alpha_i} \right) f_w.
\end{equation}  Together with $f_{w_0}=1$ this gives
\begin{equation}
f_w=\prod_{\substack{\alpha \in R^+ \\ \alpha \notin R(w)}} (1-f_\alpha),
\end{equation} where 
$$R(w)=\{\alpha \in R^+ \ | w(\alpha) \in R^- \}$$ is the inversion set of $w$.  Hence
\begin{equation}
\sum_{w \in S_n} w=\sum_{w \in S_n} \sigma_w \prod_{\alpha \notin
R(w)} (1-f_\alpha).
\end{equation} 

Now, recalling that $g_\alpha \in F$ is defined by
$f_\alpha.f_{\mu,T}=g_\alpha f_{\mu,T}$, using \eqref{sigma squared} and the fact that $\sigma_i$ is self-adjoint for $\la \cdot,\cdot \ra$ gives
\begin{equation}
\la \sigma_i f_{\mu,T}, \sigma_i f_{\mu,T} \ra=\la f_{\mu,T},\sigma_i^2 f_{\mu,T} \ra=(1-g_{\alpha_i}^2) \la f_{\mu,T},f_{\mu,T} \ra
\end{equation} and hence
\begin{equation}
\la \sigma_w f_{\mu,T}, \sigma_w f_{\mu,T} \ra=\prod_{\alpha \in R(w)}
(1-g_\alpha^2) \la f_{\mu,T},f_{\mu,T} \ra,
\end{equation} so that
\begin{align*}
\la g_S,g_S \ra&=\sum_{w \in S_n} \left(\prod_{\alpha \notin
R(w)} (1-g_\alpha)\right)^2 \la \sigma_w(f_{\mu,T}),\sigma_w(f_{\mu,T}) \ra \\
&=\sum_{w \in S_n} \left(\prod_{\alpha \notin R(w)} (1-g_\alpha)
\right)^2
\prod_{\alpha \in R(w)} (1-g_\alpha^2) \la f_{\mu,T},f_{\mu,T} \ra \\
&=n! \prod_{\alpha
\in R^+} (1-g_\alpha) \la f_{\mu,T},f_{\mu,T} \ra
\end{align*}  where the last equality depends upon a specialization of
the identity 
$$\sum_{w \in S_n} \prod_{\epsilon_i-\epsilon_j \in R(w)} \left(
1+\frac{r c_0}{z_i-z_j} \right) \prod_{\epsilon_i-\epsilon_j
\in R^+-R(w)} \left( 1-\frac{rc_0}{z_i-z_j} \right) =n!$$  This last identity follows by
(1) observing that the left-hand side is invariant by each simple
reflection $s_i$, (2) observing that multiplying it by the
discriminant clears its denominator, and (3) using the fact that the
discriminant is the minimal degree alternating polynomial.  In the formula for $\la f_{\mu,T},f_{\mu,T} \ra$ from
part (c) of Theorem~\ref{Upper triangular} we reindex by setting $b=T^{-1} w_\mu(i)$ so that $\mu_i=S(b)$, and take products over $b \in \lambda^\bullet$ instead of $1 \leq i \leq n$.  Combined with the preceding calculation,  a bit of rearranging finishes the computation of $\la g_S,g_S \ra$.
\end{proof}

\section{The minimal degree symmetric polynomial and aspherical values}

Throughout this section we fix an $r$-partition $\lambda^\bullet$ of
$n$, and write $S$ for the column-strict tableau on $\lambda^\bullet$ with
$S(b)=l+(i-1)r$ if $b$ is in the $i$th row of $\lambda^l$.  The
function $g_S$ is the (unique up to scalars)
minimal degree element of $M(\lambda^\bullet)^W$, and we set
\begin{equation}
g_{\lambda^\bullet}=g_S
\end{equation} with $S$ as above.  We have $g_S=e.f_{\mu,T}$ for a Young tableau $T$ that is uniquely determined by $S$, and where $\mu$ is the non-decreasing rearrangement of the entries of $S$.  As in the proof of Theorem~\ref{norm formula}, we can construct $g_S=e.f_{\mu,T}$ by rewriting $e$ as a certain $F$-linear combination of intertwining operators $\sigma_w$, and it follows that the leading term of $g_S$ is $n_\mu x^{\mu^+} w_0 v_T$, where $n_\mu$ is the size of the stabilizer of $\mu$ in $S_n$.  On the other hand, since the space of $W$-invariants in this degree is $1$-dimensional by Theorem~\ref{norm formula}, $g_{\lambda^\bullet}$ must be a multiple of its specialization to $c_0=d_0=\cdots=d_{r-1}=0$ (which is well-defined since the intertwiners $\sigma_i$ have no poles when $c_0=0$).  Since they have the same leading term they must be equal.  Thus while any non-symmetric $f_{\mu,T}$ with $e.f_{\mu,T}=g_{\lambda^\bullet}$ may depend on the parameters, its symmetrization does not.  

\subsection{The norm of the minimal degree symmetric function}

Let $\mathrm{lrim}^l$ be the \emph{lower rim} of $\lambda^l$, consisting
of those boxes $b \in \lambda^l$ that are not directly above another
box of $\lambda^l$.  Similarly, let $\mathrm{rrim}^l$ be the \emph{right rim} of
$\lambda^l$, consisting of those boxes which are not directly to the
left of another box of $\lambda^l$.  The \emph{lower rim}
$\mathrm{lrim}(\lambda^\bullet)$ and \emph{right rim}
$\mathrm{rrim}(\lambda^\bullet)$ of $\lambda^\bullet$ are the unions of
the respective rims of its components.  The \emph{hook length product}
associated to $\lambda^\bullet$ is
\begin{equation}
H_{\lambda^\bullet}=\prod_{\substack{b \in
\mathrm{lrim}(\lambda^\bullet) \\ b' \in \mathrm{rrim}(\lambda^\bullet)}}
\quad \prod_{\substack{1 \leq k \leq S(b)-S(b') \\ k=\beta(b)-\beta(b') \
\mathrm{mod} \ r}} (k-(d_{\beta(b)}-d_{\beta(b')})-r
(\mathrm{ct}(b)-\mathrm{ct}(b')-1)c_0).
\end{equation}  In the case $r=1$ this reduces (after setting
$c_0=-\kappa$) to the hook length
product from \cite{Dun}.  For each $0 \leq l \leq r-1$ let $b_{LL}^l$ be the lower left-hand corner box of
$\lambda^l$, with the convention that if $\lambda^l=\emptyset$ then $b_{LL}^l$ lies in the $0$th row and $1$st column, and put $S_l=S(b_{LL}^l)$ (so if $\lambda^l=\emptyset$ then $S_l=l-r$) and $c_l=\mathrm{ct}(b_{LL})$.  
Define an extra product $E_{\lambda^\bullet}$ by
\begin{equation}
E_{\lambda^\bullet}=\prod_{\substack{b \in \lambda^\bullet \\ 0 \leq l
\leq r-1}} \quad \prod_{\substack{1 \leq k \leq S(b)-S_l-r \\ k=\beta(b)-l \
\mathrm{mod} \ r}}
(k-(d_{\beta(b)}-d_{l})-r (\mathrm{ct}(b)-c_l+1) c_0) 
\end{equation}
\begin{theorem} \label{minimal norm}
The norm of $g_{\lambda^\bullet}$ is given by
$$\la g_{\lambda^\bullet},g_{\lambda^\bullet} \ra=n! H_{\lambda^\bullet} E_{\lambda^\bullet}$$
\end{theorem}

\begin{proof}
As in \cite{Dun} we induct on the number of boxes in
$\lambda^\bullet$.  The base case with one box follows from
Theorem~\ref{norm formula}.  We sketch the inductive step, which is
similar to that in \cite{Dun}.  In general, suppose $\lambda^\bullet$ is
obtained from $\chi^\bullet$ by adding a box $b$ with the property
that $S(b) \geq S(b')$ for all $b' \in \chi^\bullet$.  Then $b$ is in particular a removable box, and
Theorem~\ref{norm formula} shows that
\begin{align} \label{recurrence}
\la g_{\lambda^\bullet},g_{\lambda^\bullet} \ra&=n \la
g_{\chi^\bullet},g_{\chi^\bullet} \ra \prod_{1 \leq k \leq
S(b)} (k-(d_{\beta(b)}-d_{\beta(b)-k})-r \mathrm{ct}(b) c_0) \notag \\
&\times \prod_{b' \in \chi^\bullet} \prod_{\substack{1 \leq k \leq
S(b)-S(b') \\ k=\beta(b)-\beta(b') \ \mathrm{mod} \ r }}
\frac{k-(d_{\beta(b)}-d_{\beta(b')})-r(\mathrm{ct}(b)-\mathrm{ct}(b')-1) c_0
}{k-(d_{\beta(b)}-d_{\beta(b')})-r(\mathrm{ct}(b)-\mathrm{ct}(b'))c_0} \\
&\times \prod_{b' \in \chi^\bullet} \prod_{\substack{1 \leq k \leq
S(b)-S(b')-r \\ k=\beta(b)-\beta(b') \ \mathrm{mod} \ r }}
\frac{k-(d_{\beta(b)}-d_{\beta(b')})-r(\mathrm{ct}(b)-\mathrm{ct}(b')+1) c_0
}{k-(d_{\beta(b)}-d_{\beta(b')})-r(\mathrm{ct}(b)-\mathrm{ct}(b'))c_0}. \notag
\end{align}  Fix a positive integer $k$ and define $0 \leq l \leq r-1$
by $k=\beta(b)-l\ \mathrm{mod} \ r$.  For each positive integer $i$ with $k \leq S(b)-(l+(i-1)r)$, the $i$th row of $\chi^l$
contributes a telescoping product totaling
\begin{equation} \label{- factor}
\frac{k-(d_{\beta(b)}-d_l)-r(\mathrm{ct}(b)-\mathrm{ct}(b')-1)c_0}{k-(d_{\beta(b)}-d_l)-r(\mathrm{ct}(b)+i-1) c_0}
\end{equation} to the second line in \eqref{recurrence}, where $b'$ is the rightmost box in the $i$th row.
Likewise, for each positive integer $i$ with $ k \leq
S(b)-(l+(i-1)r)-r$, the
$i$th row of $\chi^l$ contributes a
telescoping product totaling
\begin{equation} \label{+ factor}
\frac{k-(d_{\beta(b)}-d_l)-r (\mathrm{ct}(b)+i)c_0}{k-(d_{\beta(b)}-d_l)-r(\mathrm{ct}(b)-\mathrm{ct}(b'))c_0}
\end{equation} to the third line in \eqref{recurrence} where again
$b'$ is the rightmost box in the row.  For $1 \leq i \leq \frac{1}{r}(S(b)-l-k)$ the
denominator of \eqref{- factor} for the $i+1$th row and the numerator
of \eqref{+ factor} for the $i$th row cancel, and the denominator
for the first row cancels with the factor $k-(d_{\beta(b)}-d_l)-r
\mathrm{ct}(b) c_0$ from the first line of \eqref{recurrence}.  The
denominator of \eqref{+ factor} cancels a corresponding factor in $H_{\chi^\bullet}$ for the box directly
above $b$ in $\chi^\bullet$ (if there is one; otherwise the factor
\eqref{+ factor} does not appear), and the
numerator of \eqref{+ factor} for the bottom row of $\chi^l$
contributes to the expression for $E_{\lambda^\bullet}$ (if this
bottom row consists of boxes $b'$ with $S(b')<S(b)$, which can be the
case only if $l \neq \beta(b) \ \mathrm{mod} \ r$).  The remaining factors in
$E_{\lambda^\bullet}$ are those of the form
$$
k-(d_{\beta(b)}-d_{\beta(b)-k})-r \mathrm{ct}(b) c_0 \quad \hbox{for $k$ with $1 \leq k \leq S(b)$ and $\lambda^{\beta(b)-k} = \emptyset$.}
$$  These are accounted for by those factors in the first
line of \eqref{recurrence} that did not cancel with a denominator from \eqref{- factor}.  These
observations together with the inductive hypothesis prove the formula.
\end{proof}

\subsection{Another formula for the norm}
Here we record
alternative expressions for $H_{\lambda^\bullet}$ and $E_{\lambda^\bullet}$ that 
are closer in spirit to the formulas of \cite{DuOp} and \cite{Dun}.  Let ${}^t \mu$ denote the conjugate partition of the partition $\mu$, and write
$$(x)_n=x(x+1)(x+2)\cdots(x+n-1) $$ for the Pochhammer symbol.  If a partition $\mu=\emptyset$ let ${}^t\mu=\left(  0\right)  $, that is,
${}^t\mu_{1}=0$.  In the following theorem, empty products are understood to equal $1$.  
\begin{theorem}
There are $a_1,a_2 \in \CC^\times$ with
\begin{align*}
H_{\lambda^{\bullet}}  & =a_1 \prod_{k=0}^{r-1}\prod_{l=0}^{k-1}\prod_{i=1}%
^{\min\left(  {}^t \lambda_{1}^{k},{}^t \lambda_{1}^{l}\right)  }\prod_{j=1}%
^{\lambda_{i}^{k}}\left(  \frac{k-l}{r}-\frac{d_{k}-d_{l}}{r}+c_{0}\left(
{}^t \lambda_{j}^{k}+\lambda_{i}^{l}-i-j+1\right)  \right)  _{{}^t\lambda_{j}%
^{k}-i+1}\\
& \times\prod_{l=k}^{r-1}\prod_{i=1}^{\min\left( {}^t\lambda_{1}^{k},{}^t\lambda
_{1}^{l}\right)  }\prod_{j=1}^{\lambda_{i}^{k}}\left(  \frac{r+k-l}{r}%
-\frac{d_{k}-d_{l}}{r}+c_{0}\left( {}^t\lambda_{j}^{k}+\lambda_{i}%
^{l}-i-j+1\right)  \right)  _{{}^t\lambda_{j}^{k}-i}.
\end{align*}  and
\begin{align*}
E_{\lambda^{\bullet}}  & =a_2 \prod_{k=0}^{r-1}\prod_{l=0}^{k-1}\prod
_{i={}^t \lambda_{1}^{l}+1}^{{}^t\lambda_{1}^{k}}\prod_{j=1}^{\lambda_{i}^{k}}\left(
\frac{k-l}{r}-\frac{d_{k}-d_{l}}{r}+c_{0}\left(  i-j-{}^t \lambda_{1}^{l}\right)
\right)  _{i-{}^t \lambda_{1}^{l}}\\
& \times\prod_{l=k}^{r-1}\prod_{i={}^t \lambda_{1}^{l}+2}^{{}^t \lambda_{1}^{k}}%
\prod_{j=1}^{\lambda_{i}^{k}}\left(  \frac{r+k-l}{r}-\frac{d_{k}-d_{l}}%
{r}+c_{0}\left(  i-j-{}^t \lambda_{1}^{l}\right)  \right)  _{i-{}^t\lambda_{j}^{l}%
-1}.
\end{align*}
\end{theorem}
\begin{proof}
The innermost product loops of our formulae can be expressed as Pochhammer
symbols. For fixed $k,l$ such that $1\leq k,l\leq r-1$ the typical product has
the form%
\[
P_{k,l}\left(  U,A\right)  =\prod_{\substack{1\leq m\leq
U\\m=k-l\operatorname{mod}r}}\left(  m-d_{k}+d_{l}-r A c_{0}\right)  .
\]
Change the loop variable, letting $m=k-l+ri$, then $i$ has integer values and
satisfies $\frac{1-k+l}{r}\leq i\leq\frac{U-k+l}{r}$. Thus the upper limit is
$u=\left[ \frac{U-k+l}{r}\right] $, where we write $[x]$ for the
largest integer at most $x$. If $0\leq l\leq k-1$ then
$-1<-1+\frac{2}{r}\leq\frac{1-k}{r}\leq\frac{1-k+l}{r}\leq0$ and the lower
limit is $0$, and $P_{k,l}$ has $u+1$ terms; the vacuous product occurs for
$u\leq-1$. If $k\leq l\leq r-1$ then $\frac{1}{r}\leq\frac{1+l-k}{r}\leq
\frac{r-k}{r}\leq1$ and the lower limit is $1$, and $P$ has $u$ terms (vacuous
product for $u\leq0$, but it is convenient to allow $u=0$ in the following
formula). In the latter case let $i=i^{\prime}-1$ then $0\leq i^{\prime}\leq
u$; . Thus%
\begin{align*}
P_{k,l}\left(  U,A\right)    & =r^{u+1}\left(  \frac{k-l}{r}-\frac{d_{k}%
-d_{l}}{r}-c_{0}A\right)  _{u+1},0\leq l<k;\\
P_{k,l}\left(  U,A\right)    & =r^{u}\left(  \frac{r+k-l}{r}-\frac{d_{k}%
-d_{l}}{r}-c_{0}A\right)  _{u},k\leq l\leq r-1.
\end{align*}

Consider the product $H_{\lambda^{\bullet}}$. The typical $b\in\mathrm{lrim}%
(\lambda^{\bullet})$ has $\beta\left(  b\right)  =k$, coordinates $\left(
{}^t\lambda_{j}^{k},j\right)  $ for $1\leq j\leq\lambda_{1}^{k}$ and $S\left(
b\right)  =k+r\left(  {}^t\lambda_{j}^{k}-1\right)  ,\mathrm{ct}\left(  b\right)
=j-{}^t\lambda_{j}^{k}$. The typical $b^{\prime}\in\mathrm{rrim}(\lambda
^{\bullet})$ has $\beta\left(  b\right)  =l$, coordinates $\left(
i,\lambda_{i}^{l}\right)  $ for $1\leq i\leq{}^t\lambda_{1}^{l}$ and $S\left(
b^{\prime}\right)  =l+r\left(  i-1\right)  ,\mathrm{ct}\left(  b^{\prime
}\right)  =\lambda_{i}^{l}-i$. The corresponding term in $H_{\lambda^{\bullet}}$
is $P_{k,l}\left(  k-l+r\left(  {}^t\lambda_{j}^{k}-i\right)  ,j-{}^t\lambda_{j}%
^{k}+i-\lambda_{i}^{l}-1\right)  $. Evaluate $u=\left[ \frac{U-k+l}%
{r}\right] ={}^t\lambda_{j}^{k}-i$. The condition ${}^t\lambda_{j}^{k}-i\geq0$
is equivalent to $\lambda_{i}^{k}-j\geq0$. This imposes the additional bound
$i\leq{}^t\lambda_{1}^{k},$and thus the pertinent index values are $1\leq
i\leq\min\left(  {}^t\lambda_{1}^{k},{}^t\lambda_{1}^{l}\right)  ,1\leq j\leq
{}^t\lambda_{i}^{k}$. These arguments deduce the formula for $H_{\lambda^{\bullet}%
}$ from the one in the Theorem.

Consider the product $E_{\lambda^{\bullet}}$. The typical $b\in\lambda^{\bullet}$
has~$\beta\left(  b\right)  =k,$ coordinates $\left(  i,j\right)  $ with
$1\leq i\leq{}^t\lambda_{1}^{k},1\leq j\leq\lambda_{i}^{k}$ and $S\left(
b\right)  =k+r\left(  i-1\right)  ,\mathrm{ct}\left(  b\right)  =j-i$. The corresponding $b_{LL}$ has $\beta\left(
b_{LL}\right)  =l,S\left(  b_{LL}\right)  =l+r\left(  {}^t\lambda_{1}%
^{l}-1\right)  ,\mathrm{ct}\left(  b_{LL}\right)  =1-{}^t\lambda_{1}^{l}$, so
that $U=k-l+r\left(  i-{}^t\lambda_{1}^{l}\right)  -1$ and $u=i-{}^t\lambda_{1}%
^{l}-1$ and $A=j-i+{}^t\lambda_{1}^{l}$. Each pair $\left(  i,j\right)  $ with
${}^t\lambda_{1}^{l}+1\leq i\leq{}^t\lambda_{1}^{k},1\leq j\leq\lambda_{i}^{k}$
contributes to $E_{\lambda^{\bullet}}$ (when $k\leq l$ and $i={}^t\lambda_{1}%
^{l}+1$ the contribution to the product is $1)$.  This demonstrates
the alternative formula for $E_{\lambda^\bullet}$.  
\end{proof}

\subsection{Aspherical values}

In this subsection we will, for the first time, specialize parameters to complex numbers.

A parameter $c=(c_0,d_0,\dots,d_{r-1})$ is \emph{aspherical} if $\HH_c
e_+ \HH_c \neq \HH_c$, where $e_+=\frac{1}{|W|}\sum_{w \in W} w$ is
the trivial idempotent for $W=G(r,1,n)$.  By Theorem 4.1 of \cite{BeEt}, $c$ is
aspherical exactly if there is an $r$-partition $\lambda^\bullet$ of
$n$ so that $L_c(\lambda^\bullet)^W=0$. 

We thank Pavel Etingof for pointing out the following
lemma and that it follows from S. Montarani's work \cite{Mon}.  We
give a proof based on \cite{Gri} together with \cite{BeEt}.
\begin{lemma} \label{aspherical lemma}
\begin{enumerate}
\item[(a)]  If there is a rectangle $\lambda$ with at most $n$ boxes, an integer $0 \leq l
\leq r-1$, and a positive integer $k$ with
$$
k=d_l-d_{l-k}+r \mathrm{ct}(b) c_0, \ k \neq 0 \ \mathrm{mod} \ r, \quad
\mathrm{and} \quad 1 \leq k \leq l+(\mathrm{row}(b)-1)r
$$ where $b$ is the removable (lower right hand corner) box of $\lambda$, then $c$ is aspherical.
\item[(b)] Parameters $c$ such that $c_0=-k/m$ for integers $1
\leq k< m \leq n$ are aspherical.
\end{enumerate}
\end{lemma} 
\begin{proof}
(a) Consider the module $M_c(\lambda^\bullet)$ for the rational Cherednik algebra
of type $G(r,1,m)$, where $m$ is the number of boxes in the rectangle
and $\lambda^\bullet$ is the $r$-partition with $\lambda^l=\lambda$
and all other components empty.  Provided that $c_0$ is not a rational
number, Theorem 7.5 of \cite{Gri} implies that the span of the
$f_{\mu,T}$ for $T \in \mathrm{SYT}(\lambda^\bullet)$ and $\mu_i \geq k$
for some $1 \leq i \leq m$ is a submodule of $M_c(\lambda^\bullet)$;
in particular, by the
inequality $k \leq l+(\mathrm{row}(b)-1) r$ all the polynomials $g_S$ are in this submodule.  Therefore $\la g,h \ra=0$
for all symmetric functions $g,h \in M_c(\lambda^\bullet)^W$.  By
continuity this condition continues to hold for all parameters $c$ on
the hyperplane, implying
that such $c$ are aspherical for the rational Cherednik algebra for
$G(r,1,m)$.  By Theorem 4.1 of \cite{BeEt}, it is aspherical also for
the rational Cherednik algebra of type $G(r,1,n)$.  

(b) These points are (some of the) zeros of the norm of $g_{\lambda^\bullet}$, where
$\lambda^{r-1}=(1,1,\dots,1)$ and $\lambda^l=\emptyset$ for $l \neq
r-1$ corresponds to the determinant representation of
$G(r,1,n)$.  So they are aspherical because every invariant in
$M(\lambda^\bullet)$ is a multiple of $g_{\lambda^\bullet}$.
\end{proof}

In other words, in addition to the hyperplanes $c_0=-k/m$ with $1 \leq
k<m \leq n$, for each integer $m$ with $-(n-1) \leq m \leq n-1$ the hyperplanes
\begin{equation} \label{aspherical hyperplanes}
k=d_l-d_{l-k}+rmc_0 \quad \mathrm{with} \quad 1 \leq k \leq l+\left(\sqrt{n+\frac{1}{4} m^2}-\frac{1}{2} m-1\right)r
\end{equation} are also aspherical.

\begin{theorem} \label{aspherical theorem}
The set of aspherical values $c$ for the rational Cherednik algebra
for $G(r,1,n)$ is the union of the hyperplanes $c_0=-k/m$ for integers
$1 \leq k<m \leq n$ and hyperplanes
$$k=d_l-d_{l-k}+r \mathrm{ct}(b) c_0$$ where $b$ is the lower right-hand
corner box of a rectangle with at most $n$ boxes, $0 \leq l \leq r-1$,
$k \neq 0\ \mathrm{mod} \ r$, and $1 \leq k \leq l+(\mathrm{row}(b)-1)r$.
\end{theorem}
\begin{proof}
By Lemma~\ref{aspherical lemma} the stated hyperplanes consist of
aspherical values.  Conversely, assume that
$c=(c_0,d_0,\dots,d_{r-1})$ is aspherical and choose an $r$-partition
$\lambda^\bullet$ of $n$ so that $L_c(\lambda^\bullet)^W=0$.  

Writing $g_{\lambda^\bullet}=\sum g_T v_T$ for certain $g_T \in \CC[x_1,\dots,x_n]$ we have
$$\la g_{\lambda^\bullet}, g_{\lambda^\bullet} \ra=\sum \la v_T,
\overline{g_T}.g_{\lambda^\bullet} \ra,$$ where $\overline{\cdot}$ is
the $W$-equivariant conjugate linear isomorphism of
$\CC[x_1,\dots,x_n]$ onto $\CC[y_1,\dots,y_n]$ determined by
$\overline{x_i}=y_i$.  While the left hand side is only defined for real values of the parameters (as otherwise the contravariant form does not specialize appropriately), the right hand side always makes sense and is a polynomial in the parameters $(c_0,d_0,\dots,d_{r-1})$.  For generic real  values of the parameters it is given by the formula in Theorem~\ref{minimal norm}.  It follows that it is given by this formula for all values of the parameters.   

Thus four situations can occur, corresponding to the zeros of the formula in Theorem~\ref{minimal norm}: 
\begin{enumerate}
\item[(a)] If one of the factors of $H_{\lambda^\bullet}$ corresponding to boxes $b,b'$ with $\beta(b)=\beta(b')$ is zero, then number $c_0$ is rational, of the form $c_0=-k/m$ with
integers $1 \leq k < m \leq n$.

\item[(b)]  If one of the factors in $E_{\lambda^\bullet}$
corresponding to an empty $\lambda^j$ gives a zero, then there is a box $b \in \lambda^\bullet$ and an integer $1 \leq k \leq S(b)$ (necessarily $\beta(b)=S(b) \neq j \ \mathrm{mod} \ r$) with $k=\beta(b)-j \ \mathrm{mod} \ r$ such that
$$k=d_{\beta(b)}-d_{\beta(b)-k}+r \mathrm{ct}(b) c_0,$$ in which case
putting $l=\beta(b)$ and taking the rectangle with corner at $b$ does the trick.  

\item[(c)]  There are boxes $b \in \mathrm{lrim}(\lambda^\bullet)$ and
$b' \in \mathrm{rrim}(\lambda^\bullet)$ and an integer $k$ with $1 \leq
k \leq S(b)-S(b')$ and $0 \neq k=\beta(b)-\beta(b')\ \mathrm{mod} \ r$ so that
$$k=d_{\beta(b)}-d_{\beta(b')}-r(\mathrm{ct}(b)-\mathrm{ct}(b')-1)c_0.$$  Let $l=\beta(b)$ and let $m=\mathrm{ct}(b)-\mathrm{ct}(b')-1$.  Using \eqref{aspherical hyperplanes} it will suffice to show that 
$$S(b)-S(b') \leq l+\left(\sqrt{n+\frac{m^2}{4}}-\frac{m}{2}-1\right)r.$$  Write 
$$x=\mathrm{row}(b), \ y=\mathrm{col}(b), \ x'=\mathrm{row}(b'), \ \mathrm{and} \ y'=\mathrm{col}(b').$$  Since $S(b)-S(b')=l-\beta(b')+(x-x')r$ it suffices to show that
$$x-x' \leq
\sqrt{n+\frac{(y-x-y'+x'-1)^2}{4}}-\frac{y-x-y'+x'-1}{2}-1,$$ or,
after rearranging, that
$$(x-x'+1)(y-y') \leq n .$$  This inequality follows from $xy \leq n$ and $x'y' \leq n$.  
 
\item[(d)]  There is a box $b \in \lambda^\bullet$, an integer $0 \leq
i \leq r-1$ so that $\lambda^i \neq \emptyset$ and if $b'$ is the
lower left hand corner of $\lambda^i$ then $S(b')<S(b)$, and
an integer $1 \leq k \leq S(b)-S(b')-r$ with $0 \neq
k=\beta(b)-\beta(b') \ \mathrm{mod} \ r$
and
$$k=d_{\beta(b)}-d_{\beta(b')}+r (\mathrm{ct}(b)-\mathrm{ct}(b')+1) c_0.$$
Setting
$m=\mathrm{ct}(b)-\mathrm{ct}(b')+1$ it suffices to show that
$$S(b)-S(b')-r \leq l +
\left(\sqrt{n+\frac{m^2}{4}}-\frac{m}{2}-1\right)r,$$ and the rest of the proof of this case proceeds as in case (c).
\end{enumerate}
\end{proof}  

\subsection{Other linear characters}

Let $W$ be a complex reflection group, and with the notation of section~5 of Rouquier's paper \cite{Rou} let $\HH_h$ be the rational
Cherednik algebra with parameter $h$ attached to $W$.  Fix a linear
character $\xi:W \rightarrow \CC^\times$ and write $e_\xi=\sum_{w \in
W} \xi(w)^{-1} w$ for the
corresponding symmetrizer.  With $\theta_\xi$ as in section 3.3.1 of
\cite{Rou}, there is an isomorphism
$\HH_h \rightarrow \HH_{\theta_\xi(h)}$ such that $e_\xi \mapsto
e_{\mathrm{triv}}$, and it follows that $\HH_h e_\xi \HH_h \neq \HH_h$
exactly if $\HH_{\theta_\xi(h)} e_{\mathrm{triv}} \HH_{\theta_\xi(h)}
\neq \HH_{\theta_\xi(h)}$.  

In case $W=G(r,1,n)$, a character $\xi$ is
determined by the the integers $0 \leq i \leq 1$ and $0 \leq j \leq
r-1$ with $\xi(s_1)=(-1)^i$ and $\xi(\zeta_1)=\zeta^j$.  In terms of
our parametrization, the parameter $\theta_\xi(h)$ is obtained by replacing
$c_0$ by $(-1)^i c_0$ and $d_l$ by $d_{l+j}$.  We deduce the following
corollary.
\begin{corollary} \label{twists}
Let $\xi$ be the linear character of $G(r,1,n)$ determined by
$\xi(s_1)=(-1)^i$ and $\xi(\zeta_1)=\zeta^j$.  A parameter
$c$ is such that $\HH_c e_\xi \HH_c \neq \HH_c$ exactly if one of the
following holds:
\begin{enumerate}
\item[(a)] $c_0=(-1)^{i+1} k/m$ for integers $1 \leq k < m \leq n$, or
\item[(b)] there are integers $0 \leq l \leq r-1$ and $k$ and a
rectangle with at most $n$ boxes so that writing $b$ for its lower
right-hand corner box, $$1 \leq k \leq l+(\mathrm{row}(b)-1)r \quad
\mathrm{and} \quad k=d_{l+j}-d_{l+j-k}+(-1)^i r\mathrm{ct}(b)c_0.$$
\end{enumerate}
\end{corollary}

\subsection{The groups $G(r,p,n)$}

Fix a positive integer $p$ dividing $r$ and let $G(r,p,n)$ be the subgroup of $G(r,1,n)$ consisting of matrices so
that the product of the non-zero entries is an $r/p$-th root of
unity.  When $n \geq 3$ (to avoid fusion of conjugacy classes of
reflections) and $d_i=d_j$ for $i=j \ \mathrm{mod} \ r/p$, the
rational Cherednik algebra for $G(r,p,n)$ is the subalgebra of that
for $G(r,1,n)$ generated by $\CC[x_1,\dots,x_n]$, $G(r,p,n)$, and
$\CC[y_1,\dots,y_n]$.  As in section 9 of \cite{Gri}, it is the fixed
subalgebra for a cyclic group of automorphisms of $\HH_c$, which
allows one to apply (an especially simple version of) Clifford theory to relate
representations of the two algebras.

Let $C$ be the cyclic shift (by $r/p$) operator on $r$-partitions given by
$$C.(\lambda^0,\dots,\lambda^{r-1})=(\lambda^{r/p},\lambda^{r/p+1},\dots,\lambda^{r/p-1}).$$
 Thus the representations of $G(r,1,n)$ whose restriction to
$G(r,p,n)$ contains the trivial
representation are indexed by precisely those $r$-partitions in the
$C$-orbit of $((n),\emptyset,\dots,\emptyset)$. 

With notation as in section 9 of \cite{Gri}, the equation
$$L_c(\lambda^\bullet)=\bigoplus_{q=0}^{p/k-1}
L_c(\lambda^\bullet,q)$$ implies that
$L_c(\lambda^\bullet,q)^{G(r,p,n)}=0$ for some $q$ exactly if
$L_c(\lambda^\bullet)^\xi=0$ for all linear characters $\xi$ of
$G(r,1,n)$ that restrict to the trivial linear character of
$G(r,p,n)$.  

With the specialization $d_i=d_j$ when $i=j \ \mathrm{mod} \ r/p$,
Corollary~\ref{twists} implies that $c$ is $\xi$-aspherical for some $\xi$
restricting trivially to $G(r,p,n)$ exactly if it is aspherical.  Therefore
$c$ is aspherical for the rational Cherednik algebra of type
$G(r,p,n)$ exactly if it is so for $G(r,1,n)$.  We obtain:

\begin{corollary} \label{p corollary}
For $n \geq 3$, the aspherical set for the rational Cherednik algebra
of type $G(r,p,n)$ is given by the equations of
Theorem~\ref{aspherical theorem}
restricted to parameters with $d_i=d_j$ for $i=j \ \mathrm{mod} \
r/p$.  
\end{corollary} 

\section{Ordering category $\OO_c$}

\subsection{The ordering, numerically}
For a non-zero number $m \in \RR^\times$ and $a,b \in \RR$ write $a=b\ \mathrm{mod} \ m$ if
$(a-b)m^{-1} \in \ZZ$.  For multisets $X$ and $Y$ of real numbers write
$X=Y\ \mathrm{mod} \ m$ if the corresponding multisets of equivalence classes
are equal.  Let $c=(c_0,d_1,\dots,d_{r-1}) \in \RR^r$ be a parameter
with $c_0>0$ and define a partial order $\geq_c$ on
the set of $r$-partitions of $n$ by the rule: $\lambda^\bullet \geq_c \chi^\bullet$ if for all $j \in \RR$ and $0 \leq l \leq r-1$,
\begin{align} \label{dominance condition}
&|\{b \in \lambda^\bullet \ | \ \frac{d_{\beta(b)}}{r c_0}+\mathrm{ct}(b) > j \ \mathrm{or} \ \frac{d_{\beta(b)}}{r c_0}+\mathrm{ct}(b)=j \ \mathrm{and} \ \beta(b) \leq l \}| \\
& \geq |\{b' \in \chi^\bullet \ | \ \frac{d_{\beta(b')}}{r
c_0}+\mathrm{ct}(b') > j \ \mathrm{or} \ \frac{d_{\beta(b')}}{r
c_0}+\mathrm{ct}(b')=j \ \mathrm{and} \ \beta(b') \leq l \}| \notag.
\end{align}  Also define an equivalence relation $\equiv_c$ by
$\lambda^\bullet \equiv_c \chi^\bullet$ if there is an equality
of multisets
\begin{equation} \label{core condition}
\{\mathrm{ct}(b)+\frac{d_{\beta(b)}-\beta(b)}{r c_0} \ | \ b \in \lambda^\bullet \}=\{\mathrm{ct}(b)+\frac{d_{\beta(b)}-\beta(b)}{rc_0} \ | \ b \in
\chi^\bullet \} \quad \mathrm{mod} \ c_0^{-1}.
\end{equation} 

For the definition of a highest weight category, see \cite{Rou},
Definition 4.11.  The congruence condition part of the following
theorem is analogous to a result of \cite{GrLe} concerning cyclotomic
Hecke algebras.  In \cite{LyMa} Lyle and Mathas finished the classification of the blocks of cyclotomic Hecke algebras (and showed that they are the same as those for the affine Hecke algebra).  It might be that the norm formulas presented in this paper are enough to give an alternative proof of this theorem; we have not yet attempted to do so.
\begin{theorem} \label{ordering theorem}
For a parameter $c=(c_0,d_1,\dots,d_{r-1}) \in \RR^r$ with $c_0>0$ category
$\OO_c$ is a highest weight category with respect to the order given
by $\lambda^\bullet \geq \chi^\bullet$ if
$\lambda^\bullet \geq_c \chi^\bullet$ and $\lambda^\bullet \equiv_c \chi^\bullet$.    
\end{theorem}
\begin{proof}
We will show that if $[M_c(\lambda^\bullet):L_c(\chi^\bullet)] \neq 0$
then $\lambda^\bullet \geq_c \chi^\bullet$ and $\lambda^\bullet
\equiv_c \chi^\bullet$.  Assuming this for the moment, for each $r$-partition
$\lambda^\bullet$, write $P_c(\lambda^\bullet)$ for the projective cover
of $L_c(\lambda^\bullet)$.  By Corollary 2.10 of \cite{GGOR}
$M_c(\lambda^\bullet)$ is a quotient of $P_c(\lambda^\bullet)$.  It follows from the
formula following Theorem 2.19 and Proposition 3.3 of \cite{GGOR} that
if $[P_c(\chi^\bullet):M_c(\lambda^\bullet)] \neq 0$ then $\lambda^\bullet
\geq \chi^\bullet$, and hence $\OO_c$ is a highest weight category with
respect to the order $\geq$.

If $[M_c(\lambda^\bullet):L_c(\chi^\bullet)] \neq 0$ then upon examining generalized $\ttt$-eigenspaces part (a)
of Theorem~\ref{Upper triangular} implies that there exist orderings $b_1,b_2,\dots,b_n$ and $b_1',b_2',\dots,b_n'$ of the boxes of $\lambda^\bullet$ and $\chi^\bullet$, respectively, and non-negative integers $\mu_i$ so that
\begin{equation} \label{congruence condition}
\beta(b_i)-\mu_i=\beta(b_i') \ \mathrm{mod} \ r
\end{equation} and
\begin{equation} \label{domcore condition}
\mu_i=d_{\beta(b_i)}-d_{\beta(b_i')}+r(\mathrm{ct}(b_i)-\mathrm{ct}(b_i')) c_0
\end{equation} for $1 \leq i \leq n$.  The latter condition can be rewritten
\begin{equation} \label{eqn a}
0 \leq \mu_i=\left(\frac{d_{\beta(b_i)}}{r c_0}+\mathrm{ct}(b_i)-\left(\frac{d_{\beta(b_i')}}{r c_0}+\mathrm{ct}(b_i') \right) \right) r c_0,
\end{equation}  and setting $\mu_i=q_i r+\beta(b_i)-\beta(b_i')$ one obtains also
\begin{equation} \label{eqn b}
0 \leq q_i=\left(\frac{d_{\beta(b_i)}-\beta(b_i)}{r c_0}+\mathrm{ct}(b_i)-\left(\frac{d_{\beta(b_i')}-\beta(b_i')}{r c_0}+\mathrm{ct}(b_i') \right) \right) c_0
\end{equation}  This implies $\lambda^\bullet \equiv_c \chi^\bullet$,
and also since $c_0>0$, equations \eqref{eqn a} and \eqref{eqn b} imply that for $1 \leq i \leq n$
\begin{equation} \label{cond}
\frac{d_{\beta(b_i)}}{r c_0}+\mathrm{ct}(b_i) \geq \frac{d_{\beta(b_i')}}{r c_0}+\mathrm{ct}(b_i') \quad \hbox{with equality implying} \quad \beta(b_i) \leq \beta(b_i').
\end{equation}  Hence for all $j \in \RR$ and $0 \leq l \leq r-1$,
\begin{align*}
&|\{b \in \lambda^\bullet \ | \ \frac{d_{\beta(b)}}{r c_0}+\mathrm{ct}(b) > j \ \mathrm{or} \ \frac{d_{\beta(b)}}{r c_0}+\mathrm{ct}(b)=j \ \mathrm{and} \ \beta(b) \leq l \}| \\
& \geq |\{b' \in \chi^\bullet \ | \ \frac{d_{\beta(b')}}{r c_0}+\mathrm{ct}(b') > j \ \mathrm{or} \ \frac{d_{\beta(b')}}{r c_0}+\mathrm{ct}(b')=j \ \mathrm{and} \ \beta(b') \leq l \}|,
\end{align*} which is \eqref{dominance condition}.
\end{proof}

\subsection{Cores, quotients, and beta numbers} \label{cores subsection}

The basic reference we use for this material is Chapter 1 of
\cite{Mac}.  See especially exercise 8 of section~1.  

Let $s$ be a complex number and let $\lambda$ be a partition.  The set of \emph{beta numbers} of $\lambda$ with respect to $s$ is
\begin{equation}
B_s(\lambda)=\{\lambda_j+s-j+1 \ | \ 1 \leq j < \infty \}.
\end{equation}  We are using the notation $B_s$ to avoid a conflict with the function $\beta$ we defined on boxes of an $r$-partition.  The first $\ell(\lambda)$ beta numbers of $\lambda$ are $s+1$ more than the contents of the rightmost boxes of $\lambda$.  Since no two of the rightmost boxes of $\lambda$ can have the same content, $\lambda$ and $s$ can be reconstructed from the beta numbers.

The reflection representation of $S_r=W(A_{r-1})$ is
\begin{equation*}
\{(a_1,a_2,\dots,a_r) \in \RR^r \ | \ \sum_{i=1}^r a_i=0 \} \ \hbox{and the root lattice is} 
\ Q=\{(a_1,\dots,a_r) \in \ZZ^r \ | \ \sum_{i=1}^r a_i=0 \}.
\end{equation*} 

Let $\lambda^\bullet=(\lambda^0,\lambda^1,\dots,\lambda^{r-1})$ be an
$r$-partition and let $a=(a_1,a_2,\dots,a_r)$ be an element of the
root lattice $Q$ of type $A_{r-1}$.  As in Section~6 of \cite{Gor}, we will construct from this data a partition $\lambda$ whose $r$-quotient is $\lambda^\bullet$ and whose $r$-core is determined by $a$.  In order to conform with the notation in Gordon's paper we set $\lambda^{(i)}=\lambda^{r-i}$, for $1 \leq i \leq r$.  The set of $B_0$ numbers of $\lambda$ is
\begin{equation} \label{B0 numbers}
B_0(\lambda)=\bigcup_{1 \leq i \leq r} \{i+r(x-1) \ | \ x \in B_{a_i}(\lambda^{(i)}) \}.
\end{equation}  This defines a bijection
\begin{equation}
Q \times \Pi^r \rightarrow \Pi,
\end{equation} where $\Pi$ is the set of partitions.

\subsection{Dominance order}
If $\lambda$ and $\chi$ are (ordinary) partitions of $n$ then we write $\lambda \geq \chi$ if for all $1 \leq i \leq n$ we have
\begin{equation} \label{dom def}
\sum_{j=1}^i \lambda_j \geq \sum_{j=1}^i \chi_j.
\end{equation}  The relation $\geq$ is called \emph{dominance order}.  Note that we denoted it by $\geq_d$ in Section~1 in order to distinguish it from the partial order we defined on $\ZZ_{\geq 0}$, but since the latter ordering plays no role from now on we henceforth drop the subscript $d$ and reserve the symbol $\geq$ for dominance order.

The following lemma gives the relationship between dominance order and
the special case $r=1$ of Theorem~\ref{ordering theorem}; we will use it
to prove Corollary~\ref{iainquestion}. 

\begin{lemma}
Suppose $\lambda$ and $\chi$ are partitions of $n$.  Then $\lambda \geq \chi$ if and only if there exist numbering $b_1,b_2,\dots,b_n$ and $b_1',b_2',\dots,b_n'$ of the boxes of $\lambda$ and $\chi$ such that $\mathrm{ct}(b_i) \geq \mathrm{ct}(b_i')$ for $1 \leq i \leq n$.  Equivalently, $\lambda \geq \chi$ if and only if for each $j \in \ZZ$, one has
\begin{equation} \label{ct comparison}
|\{b \in \lambda \ | \ \mathrm{ct}(b) \geq j \}| \geq |\{b' \in \chi \ | \ \mathrm{ct}(b') \geq j \}|
\end{equation} 
\end{lemma}
\begin{proof}
Suppose $\lambda \geq \chi$.  By \cite{Mac} (1.15) we may obtain
$\lambda$ from $\chi$ by applying a sequence of raising operators (that
is, operators that move a single box to a higher row).  But
when $\lambda=R \chi$ for a single raising operator $R$ the equation
\eqref{ct comparison} is clear.  

Conversely, assuming that numberings as in the statement of the lemma exist, it follows that for each $c \in \ZZ$ there are at least as many boxes of $\lambda$ with content at least $c$ as boxes of $\chi$ with content at least $c$.  For $1 \leq i \leq n$ such that $\lambda_i>0$, let $c_i=\lambda_i-i$ be the content of the last box in the $i$th row of $\lambda$ and compute
\begin{align*}
\sum_{j=1}^i (\lambda_j-\chi_j)&=|\{\hbox{boxes in first $i$ rows of $\lambda$ with content at least $c_i$} \}| \\
&-|\{\hbox{boxes in first $i$ rows of $\chi$ with content at least $c_i$} \}| \\
&+|\{\hbox{boxes in first $i$ rows of $\lambda$ with content less than $c_i$} \}| \\
&-|\{\hbox{boxes in first $i$ rows of $\chi$ with content less than $c_i$} \}| \\
&\geq |\{\hbox{boxes in $\lambda$ with content at least $c_i$} \}| \\
&-|\{\hbox{boxes in $\chi$ with content at least $c_i$} \}| \geq 0.
\end{align*}  Here we used the fact that 
\begin{equation*}
\{\hbox{boxes in first $i$ rows of $\lambda$ with content at least $c_i$} \}=\{\hbox{boxes in $\lambda$ with content at least $c_i$} \}.
\end{equation*}
\end{proof}
If $\lambda$ is a partition then for $k \in \ZZ$ and $s \in \RR$, the relationship between beta numbers and contents described in \ref{cores subsection} implies
\begin{equation} \label{beta and contents}
|\{b \in \lambda \ | \ \mathrm{ct}(b) = k \}|=\begin{cases} |\{x \in B_s(\lambda) \ | \ x \geq k+s+1 \}| \quad &\hbox{if $k \geq 0$, and} \\
 |\{x \in B_s(\lambda) \ | \ x \geq k+s+1 \}| +k \quad &\hbox{if $k<0$.}
\end{cases}
\end{equation}  We set $a_i=d_{r-i}/r c_0$ and assume for the
remainder of this section that $a_i \in \ZZ$ for $1 \leq i \leq r$.  If $\lambda$ is the partition corresponding to $\lambda^\bullet$ and $a=(a_1,\dots,a_r)$ under the bijection of Section~\ref{cores subsection}, then for $j \in \RR$ using \eqref{beta and contents} gives
\begin{align*}
|\{b \in \lambda \ | \ \mathrm{ct}(b) \geq j \}|
&=\sum_{\substack{ k \in \ZZ \\ k \geq j}}  \sum_{1 \leq l \leq r} |\{x \in B_{a_l}(\lambda^{(l)}) \ | \ x \geq (k-l+1)/r+1 \}| + \sum_{k \in \ZZ, \ j \leq k <0} k \\ 
&=\sum_{\substack{ k \in \ZZ \\ k \geq j}} 
\sum_{1 \leq l \leq m_k+1} |\{x \in B_{a_l}(\lambda^{(l)}) \ | \ x \geq (k-m_k)/r +1 \}| \\
&+\sum_{m_k +1 < l \leq r} |\{x \in B_{a_l}(\lambda^{(l)}) \ | \ x \geq (k-m_k)/r \}|  
+\sum_{k \in \ZZ, \ j \leq k <0} k \\ 
\end{align*} where for $k \in \ZZ$ we write $m_k$ for the integer with $0 \leq m_k \leq r-1$ and $k=m_k \ \mathrm{mod} \ r$.  Therefore using \eqref{beta and contents} again
\begin{align*}
|\{b \in \lambda \ | \ \mathrm{ct}(b) \geq j \}|
&=\sum_{\substack{ k \in \ZZ \\ k \geq j}} 
\sum_{1 \leq l \leq m_k+1} |\{b \in \lambda^{(l)} \ | \ \mathrm{ct}(b)=(k-m_k)/r -a_l \}| \\
&+\sum_{m_k+1 < l \leq r} |\{b \in \lambda^{(l)} \ | \ \mathrm{ct}(b)=(k-m_k)/r -a_l-1 \}| 
 +\sum_{k \in \ZZ, \ j \leq k <0} k \\
&-\sum_{\substack{k \in \ZZ, k \geq j \\ 1 \leq l \leq m_k+1 \\ (k-m_k)/r-a_l<0}} \left( (k-m_k)/r -a_l \right) -\sum_{\substack{k \in \ZZ, k \geq j \\ m_k+1 < l \leq r \\
     (k-m_k)/r -a_l-1<0}} \left( (k-m_k)/r-a_l-1
 \right) 
\end{align*}  Note that the last three summands depend on the sequence
$a_1,\dots,a_r$ and the number $j$ but are independent of
$\lambda$.  Write $f(a_\bullet,j)$ for their sum and $n_j=r q_j+m_j$
with $n_j$ the least integer greater than or equal to $j$ and
$q_j,m_j$ integers with $0 \leq m_j \leq r-1$.  In the sums over $k$ above we reindex by setting $q=(k-m_k)/r$ to obtain
\begin{align*}
|\{b \in &\lambda \ | \ \mathrm{ct}(b) \geq j \}|-f(a_\bullet,j)\\
&=\sum_{m=m_j+1}^{r} \sum_{q \geq q_j} \left( \sum_{l=1}^m |\{b \in
\lambda^{(l)} \ | \ \mathrm{ct}(b)=q-a_l \}|+\sum_{l=m+1}^r |\{b \in
\lambda^{(l)} \ | \ \mathrm{ct}(b)=q-a_l-1 \}| \right)  \\
&+\sum_{m=1}^{m_j} \sum_{q \geq q_j+1} \left( \sum_{l=1}^m |\{b \in
\lambda^{(l)} \ | \ \mathrm{ct}(b)=q-a_l \}|+\sum_{l=m+1}^r |\{b \in
\lambda^{(l)} \ | \ \mathrm{ct}(b)=q-a_l -1\}| \right) \\
&=\sum_{m=1}^r \sum_{q \geq q_j+1} \left( \sum_{l=1}^m |\{b \in
\lambda^{(l)} \ | \ \mathrm{ct}(b)=q-a_l\}|+\sum_{l=m+1}^r |\{b \in
\lambda^{(l)} \ | \ \mathrm{ct}(b)=q-a_l-1 \}| \right) \\
&+\sum_{m=m_j+1}^{r} \left( \sum_{l=1}^m |\{b \in
\lambda^{(l)} \ | \ \mathrm{ct}(b)=q_j-a_l \}|+\sum_{l=m+1}^r |\{b \in
\lambda^{(l)} \ | \ \mathrm{ct}(b)=q_j-a_l-1 \}| \right).
\end{align*}  Counting the number of appearances of each summand of the form $|\{b \in \lambda^{(l)} \ | \ \mathrm{ct}(b)=X \}|$ for each integer $X \geq q_j-1-a_l$ gives
\begin{align*}
|\{b \in &\lambda \ | \ \mathrm{ct}(b) \geq j \}|-f(a_\bullet,j)\\
&=r\sum_{l=1}^r |\{b \in \lambda^{(l)} \ | \ \mathrm{ct}(b) \geq q_j+1-a_l \}| 
+\sum_{l=1}^r (l-1) |\{b \in \lambda^{(l)} \ | \ \mathrm{ct}(b)=q_j-a_l\}|  \\
&+\sum_{1 \leq l \leq m_j+1} (r-m_j) |\{b \in \lambda^{(l)} \ | \ \mathrm{ct}(b)=q_j-a_l\}|
+\sum_{m_j+1<l \leq r} (r-l+1) |\{b \in \lambda^{(l)} \ | \ \mathrm{ct}(b)=q_j-a_l\}| \\
&+\sum_{m_j+1<l \leq r} (l-m_j-1) |\{b \in \lambda^{(l)} \ | \ \mathrm{ct}(b)=q_j-1-a_l \}| \\
&=r\sum_{l=1}^r |\{b \in \lambda^{(l)} \ | \ \mathrm{ct}(b) +a_l \geq q_j+1\}| 
+\sum_{1 \leq l \leq m_j+1} (r-m_j-1+l) |\{b \in \lambda^{(l)} \ | \ \mathrm{ct}(b)+a_l=q_j \}| \\
&+r \sum_{m_j+1< l \leq r} |\{b \in \lambda^{(l)} \ | \ \mathrm{ct}(b)+a_l=q_j \}|
+\sum_{m_j+1< l \leq r} (l-m_j-1) |\{b \in \lambda^{(l)} \ | \ \mathrm{ct}(b)+a_l=q_j-1\}|
\end{align*} which, recalling $\lambda^l=\lambda^{(r-l)}$ and $a_l=d_{r-l}/rc_0$ and reindexing by interchanging $l$ and $r-l$, is equal to
\begin{align*} 
&=r \sum_{l=0}^{r-1} |\{ b \in \lambda^l \ | \ \mathrm{ct}(b)+d_l/rc_0 \geq q_j+1 \}|
+\sum_{r-m_j-1 \leq l \leq r-1} (2r-m_j-1-l) |\{b \in \lambda^l \ | \ \mathrm{ct}(b)+d_l/rc_0=q_j \}| \\
&+r \sum_{0 \leq l < r-m_j-1} |\{b \in \lambda^l \ | \ \mathrm{ct}(b)+d_l/rc_0=q_j \}| \\
&+\sum_{0 \leq l <r-m_j-1} (r-m_j-1-l) |\{b \in \lambda^l \ | \ \mathrm{ct}(b)+d_l/rc_0=q_j-1 \}| \\
&=\sum_{0\leq l <r-m_j-1} |\{b \in \lambda^\bullet \ | \ \mathrm{ct}(b) +d_{\beta(b)}/rc_0 > q_j-1, \ \text{or} \ =q_j-1 \ \text{and} \ \beta(b) \leq l \}| \\
&+\sum_{r-m_j-1 \leq l \leq r-1} |\{b \in \lambda^\bullet \ | \ \mathrm{ct}(b)+d_{\beta(b)}/rc_0  > q_j, \ \text{or} \ =q_j \ \text{and} \ \beta(b) \leq l \}|.
\end{align*} The final equality above is obtained by counting the number of times a given box $b \in \lambda^\bullet$ contributes, according to the cases $\mathrm{ct}(b)+d_{\beta(b)}/rc_0=q_j-1$, $\mathrm{ct}(b)+d_{\beta(b)}/rc_0=q_j$ and $\beta(b) < r-m_j-1$ or $\beta(b) \geq r-m_j-1$, and 
$\mathrm{ct}(b)+d_{\beta(b)}/rc_0 \geq q_j+1$.  

Write $\lambda^\bullet \geq_c' \chi^\bullet$ if $\lambda \geq \chi$ in
dominance order; note that $\lambda$ and $\chi$ depend on $c$ via the
bijection described above.  The preceding
calculation and Theorem~\ref{ordering theorem} give the answer to Question 10.1 of Gordon's paper \cite{Gor}:
\begin{corollary} \label{iainquestion}
Assume that $c_0>0$ and $d_l/rc_0 \in \ZZ$ for $0 \leq l \leq
r-1$.  Category $\OO_c$ for $\HH_c$ is a highest weight category with respect to the order $\geq_c'$.
\end{corollary}

\def\cprime{$'$} \def\cprime{$'$}


\begin{thebibliography}{GGOR}

\bibitem[BeCh]{BeCh} Y. Berest and O. Chalykh, \emph{Quasi-invariants of complex reflection groups}, arxiv:0912.4518v1

\bibitem[BeEt]{BeEt} R. Bezrukavnikov and P. Etingof, \emph{Parabolic
induction and restriction functors for rational Cherednik algebras},
 Selecta Math. (N.S.)  14  (2009),  no. 3-4, 397--425. arxiv:0803.3639

\bibitem[BEG]{BEG} Y. Berest, P. Etingof, and V. Ginzburg,
\emph{Cherednik algebras and differential operators on
quasi-invariants}, Duke Math. J.  118  (2003),  no. 2, 279--337. arXiv:math/0111005

\bibitem[Dun]{Dun} C. Dunkl, \emph{Symmetric and antisymmetric
vector-valued Jack polynomials}, arxiv:1001.4485

\bibitem[DuOp]{DuOp} C. Dunkl and E. Opdam, \emph{Dunkl operators for
complex reflection groups},  Proc. London Math. Soc. (3)  86  (2003),
no. 1, 70--108. arXiv:math/0108185

\bibitem[GGS]{GGS} V. Ginzburg, I. Gordon, and T. Stafford, \emph{Differential operators and Cherednik algebras}, Selecta Math. (N.S.) 14 (2009), no. 3-4, 629--666.

\bibitem[GGOR]{GGOR} Ginzburg, Guay, Opdam, and Rouquier,
\emph{On the category $\OO$ for rational Cherednik algebras},  Invent. Math.  154  (2003),  no. 3, 617--651. arXiv:math/0212036

\bibitem[Gor]{Gor} I. Gordon, \emph{Quiver varieties, category O for
rational Cherednik algebras, and Hecke algebras},
Int. Math. Res. Pap. IMRP  2008,  no. 3, Art. ID rpn006, 69 pp. arxiv:math/0703150

\bibitem[GoSt]{GoSt} I. Gordon and J.T. Stafford, \emph{Rational
Cherednik algebras and Hilbert schemes}, Adv. Math. 198 (2005), no. 1,
222--274, arXiv:math/0407516

\bibitem[GrLe]{GrLe} J. Graham and G. Lehrer, \emph{Cellular
algebras}  Invent. Math.  123  (1996),  no. 1, 1--34.

\bibitem[Gri]{Gri} S. Griffeth, \emph{Orthogonal functions generalizing Jack
polynomials}, to appear in Transactions of the American Mathematical
Society, arXiv:0707.0251

\bibitem[Gri2]{Gri2} S. Griffeth, \emph{Towards a combinatorial representation theory for the rational Cherednik algebra of type $G(r,p,n)$}, to appear in Proceedings of the Edinburgh Mathematical Society, arXiv:math/0612733v3

\bibitem[LyMa]{LyMa} S. Lyle and A. Mathas, \emph{Blocks of cyclotomic Hecke algebras}, Adv. Math. 216 (2007), no. 2, 854--878, arXiv:math/0607451

\bibitem[Mac]{Mac} I. Macdonald, \emph{Symmetric functions and Hall
polynomials}, second edition, Oxford Mathematical Monographs, 1995.

\bibitem[Mac2]{Mac2} I. Macdonald, \emph{Affine Hecke algebras and
orthogonal polynomials}, Cambridge University Press, 2003.

\bibitem[Mon]{Mon} S. Montarani \emph{On some finite dimensional
representations of symplectic reflection algebras associated to wreath
products},  Comm. Algebra  35  (2007),  no. 5, 1449--1467. arXiv:math/0411286

\bibitem[Rou]{Rou} R. Rouquier, \emph{$q$-Schur algebras and complex
reflection groups},  Mosc. Math. J.  8  (2008),  no. 1, 119--158,
184. arXiv:math/0509252

\end{thebibliography}
\end{document}